\documentclass{amsart}

\usepackage{t1enc}

\usepackage{fancyhdr}  

\usepackage{tikz,ifthen}
\usetikzlibrary{calc,fadings, intersections, decorations.pathreplacing}
\usepackage{tikz-3dplot}

\usepackage{hyperref}

\hypersetup{
  colorlinks   = true,   
  urlcolor     = blue,   
  linkcolor    = blue,   
  citecolor   = red      
}

\usepackage[msc-links]{amsrefs}

\newtheorem{thm}{Theorem} 
 
\newtheorem{lem}{Lemma} 
 
\newtheorem{exa}{Shembull} 
\newtheorem{exe}{Ushtrim} 
 
\newtheorem{cor}{Rrjedhim} 
\newtheorem{prob}{Problem} 
\newtheorem{pyetje}{Pyetje}

\newcommand\sol{\bigskip \noindent\textbf{Solution:} \medskip}

\def\<{\langle }
\def\>{\rangle }

\newcommand{\cvec}[1]{\begin{bmatrix} #1    \end{bmatrix}}

\newcommand\Z{\mathbb{Z}}           
\newcommand\Q{\mathbb Q}            
\newcommand\R{\mathbb R}            
\newcommand\C{\mathbb C}            

\newcommand\D{\Delta}

\def\a{\alpha}          
\def\b{\beta}
\def\s{\sigma}

\newcommand\e{\varepsilon}

\def\l{\lambda}              

\def\v{\mathbf{v}}

\def\x{\mathbf{x}}

\def\cdo{çdo }

\def\att{at\"eher\"e dhe vet\"em at\"eher\"e kur }

\def\es{\"esht\"e } 
\def\qe{q\"e }
\def\te{t\"e }

\def\ne{n\"e }

\def\at{at\"eher\"e }

\def\Nqs{N\"e qoft\"e se }
\def\nqs{n\"e qoft\"e se }
\def\att{ at\"eher\"e dhe vet\"em at\"eher\"e kur }
\def\se{s\"e }

\def\ne{n\"e }
\def\te{t\"e }
\def\qe{q\"e }

\def\proof{\noindent \textit{Vërtetim. }}

\issueinfo{10}{1}{number 1}{2016}
\copyrightinfo{2016}{\href{http://albanian-j-math.com}{Albanian Journal of Mathematics}}
\def\issn{{\sc \textbf{ISSN: }} 1930-1235; }
\def\issueyear{\textbf{2016}}
\PII{\issn  (\issueyear)}

\pagespan{47}{80}

\begin{document}

\title{Mësimdhënia e matematikës  nëpërmjet problemeve klasike}


\author{Bedri Shaska}
\address{48619 Stoneridge Dr \\  Northville, MI 48168}
\curraddr{}
\email{bshaska@risat.org}
\thanks{The first author taught mathematics in middle school and high school in Albania during the period 1960-1995. In 1996, he was awarded the medal "Mësues i popullit" from the President of Albania for his services in teaching, the highest award awarded in teaching in Albania.}

\author{Tanush Shaska}
\address{Department of Mathematics and Statistics\\ Oakland University \\ 368 Mathematics Science Center \\
Rochester MI 48309-4479
}
\curraddr{}
\email{shaska@oakland.edu}

\dedicatory{Kushtuar ish-studentëve tanë të të gjitha moshave.}

 


{\let\thefootnote\relax\footnotetext{\textit{English title:} \textbf{Teaching of mathematics through classical problems}.}}

\maketitle


\hrule

\vspace{.2cm}

\begin{abstract}
Në këtë artikull  ne trajtojmë se si mësimdhënia e matematikës në shkollat 9-vjeçare dhe të mesme mund të përmirësohet në mënyrë të ndjeshme kur motivimi i koncepteve dhe ideve bëhet nëpërmjet problemeve klasike dhe historisë së matematikës. Kjo përmirëson intuitën e studentëve, zgjon kuriozitetin e tyre, shmang mësimin përmendësh të formulave të panevojshme, dhe i vendos konceptet në këndvështrimin e duhur historik.

\end{abstract}

\begin{abstract}
In this paper we discuss how teaching of mathematics for middle school and high school students can be improved dramatically when motivation of concepts and ideas is done through the classical problems and the history of mathematics. This improves the intuition of students, awakens their curiosity, avoids memorizing useless formulas, and put concepts in a historical prospective.

\end{abstract}

\vspace{.4cm}

\hrule

\vspace{.4cm}

\noindent \textit{Mathematics Subject Classes 2010:} 94B05; 97A30; 97C70; 97D40  


\vspace{.4cm}

\hrule

\vspace{.4cm}



\section{Hyrje}

Programet mësimore të matematikës gjithmonë kanë shkaktuar diskutime jo vetëm mbi përmbajtjen e tyre, por edhe mbi këndvështrimin pedagogjik, psikollogjik, dhe filozofik që duhet të kenë programe të tilla.  Diskutime këto që kanë ndodhur jo vetëm në vendin tonë, por në gjithë vendet e tjera. Pikat kryesore të diskutimit janë të njëjta pothuajse në gjithë vendet e botës dhe përmblidhen në përgjigjen e pyetjes që vijon:  Cili është qëllimi apo objektivi i programit të matematikës në arsimin 12-vjeçar?

Pyetja kjo që ka përgjigje të ndryshme të cilat varen kryesisht nga politikat arsimore të çdo vendi.  Pak a shumë përgjigjet përmblidhen si më poshtë:  

\begin{itemize}
\item Të pregatisë nxënës me një formim të përgjithshëm matematik për të përballuar me sukses sfidat e jetës. Këtu futet kryesisht grupi i nxënësve për të cilët njohuritë matematike nuk   përbëjnë bazën e profesionit të tyre. 

\item Të pregatisë nxënës me një formimin e nevojshëm  matematik për studentët e shkencave dhe inxhinjerive.  Këtu futet kryesisht grupi i nxënësve për të cilët njohuritë matematike do të jenë baza e profesionit të tyre.

\item Të pregatisë matematikanët e ardhshëm apo mësuesit e ardhshëm të matematikës. 
\end{itemize}

Zakonisht sisteme të ndryshme arsimore përpiqen që të realizojnë të tre këto objektiva, pavarësisht se metodat dhe rruga ndryshon sipas vendeve të ndryshme.  Duhet theksuar se tashmë edhe tek profesionet ku roli i matematikës konsiderohej minimal si shkencat shoqërore, mjekësore,juridike, roli i matematikës në dekadat e fundit \es rritur në mënyrë të konsiderueshme falë rolit të statistikës, "big data", aplikimeve të matematikës në mjekësi, në shkencat biologjike, etj. 

Nga ana pedagogjike diskutimi kryesor është se cila është rruga më e mirë për të realizuar këto objektiva?  Në këtë artikull modest, ne do të përpiqemi të sygjerojmë disa ide se si të rrisim efektivitetin e mësimdhënies në orën e matematikës.

Ideja kryesore e këtij artikulli \es që mësimdhënia e matematikës  duhet të bëhet nëpërmjet problemeve klasike që kanë zhvilluar matematikën si shkencë dhe historia e zhvillimit të matematikës duhet të jetë pjesë e mësimdhënies. Kjo i jep matematikës jetë në klasë dhe e bën të jetë diçka më shumë se një grup formulash të thata. Për më tepër kjo motivon konceptet bazë dhe u jep nxënësve një arsye më shumë përse këto koncepte janë të rendësishme dhe duhen mësuar. 

Në proçesin e studimit të literaturës së këtij artikulli ne zbuluam një sasi të konsiderueshme artikujsh të tjerë që kanë mështetur këtë ide.  Më i njohuri dhe më i përkushtuari pas kësaj ideje \es Felix Klein në \cites{MR0076336, MR0015349, MR0055397, MR0080930,  MR549187, MR2098410, MR3495524,  MR3495525}. Filozofia e Klein ka dominuar pedagogjinë e matematikës në fillimet e shekullit të XX, megjithëse mund të themi se implementimi i kësaj filozofie në masë lë shumë për të dëshiruar. 
Këto shënime janë organizuar si më poshtë: 

Në Kreun~2, ne japim një përshkrim  të shkurtër të disa problemeve dhe koncepteve ku një mësues i talentuar mund të zgjerohet më shumë në program dhe si pasojë të lerë një mbresë të thellë në formimin matematik të nxënësve.  Ne nuk pretendojmë se lista e problemeve që ne sygjerojmë \es e plotë,  Një listë e tillë mund të plotësohet më tepër duket pasur parasysh  eksperiencat individuale të secilit mësues  dhe grupin e nxënësve.

Në Kreun~3 ne diskutojmë rolin kritik dhe konstruktiv të historisë së shkencës në formimin e studentëve të matematikës.  Ky rol duhet të jetë i ingranuar me programin.  Ne duhet të përpiqemi t'u shpjegojmë nxënësve jo vetëm teoritë që i mbijetuan kohës, por edhe përpjekjet e dështuara pasi nga këto dështime në shumicën e rasteve mësohet shumë. 

Në Kreun~4 ne përqëndrohemi tek problemi më klasik dhe më themelor i matematikës, ai i zgjidhjes së ekuacioneve polinomiale. Pjesa kryesore e matematikës mund të lidhet me këtë temë. Pikërisht algjebra, gjeometria, gjeometria algjebrike, një pjesë e mirë e analizës, kriptografia, teoria e kodeve, etj i kanë themelet e tyre tek polinomet ose zgjidhja e ekuacioneve polinomiale ose siç quhen ndryshe ekuacione algjebrike.  Si mund të trajtohet kjo temë në kontekstin e duhur historik duke i mësuar nxënësit jo thjesht formulën e zgjidhjes së ekuacionit kuadratik, por edhe diskriminantin, polinomet simetrike, teorinë Galua, format binare, algjebrën lineare, teorinë e invariantëve, seritë e Taylor, varietetet algjebrike, kurbat eliptike dhe hipereliptike e shumë koncepte të tjera të matematikës bashkëkohore? 

Në Kreun~5 ne trajtojmë konceptin e diskriminantit nga një këndvështrim intuitiv dhe historik. Eshtë ndoshta shembulli më i mirë që një koncept kaq intuitiv \es kthyer nga programet mësimore ne një koncept të thatë dhe në një formulë që shumica e nxënësve dhe mësuesve nuk ja kuptojnë vlerën.  Duhet pranuar se ky nuk \es vetëm problem i shkollave tona. Shumica e studentëve të matematikës mësojnë për diskriminantin në algjebrën e lartë dhe edhe atëherë thjesht formulat bazë dhe asgjë më shumë. Cila \es lidhja e diskriminantit të një polinomi dhe determinantit (ose përcaktorit) të një matrice, kur filluan këto dy koncepte të ndaheshin nga njeri-tjetri në historinë e matematikës?

Format binare trajtohen në Kreun~6. Koncepti i formave është një nga konceptet bazë të matematikës që ka luajtuar dhe vazhdon të luajë një rol të rëndësishëm në matematikë, por që çuditërisht trajtohet shumë pak në programet e matematikës së shkollës së mesme. Ne trajtojmë vetëm format kuadratike megjithëse me fare pak përpjekje gjithë teoria mund të trajtohet për gjithë format me gradë $n\geq 2$. Edhe vetëm me format kuadratike ne mund të ngremë një teori të tërë dhe të drejtojmë nxënësin tek koncepte fondamentale të matematikës si hapësirat vektoriale, veprimi i një grupi mbi një bashkësi,  ndryshimi i koordinatave, teoria e invariantëve dhe puna kolosale e matematikanëve te shekullit të XIX si Clebsch, Gordan, Bolza, Hilbert, etj. 

Format kuadratike mund të hapin edhe horizonte të tjera si teorinë e reduktimit (shih Beshaj \cite{MR3579531}), bashkësitë Julia \cite{MR3532882}, gjeometrinë hiperbolike, etj. Ne zhvillojmë aq teori në Kreun~6 sa të klasifikojmë gjithë prerjet konike dhe sipërfaqet quadratike. Që kjo të realizohet saktë ne duhet të prezantojmë nxënësit me matricat, eigenvlerat dhe eigenvektorët, procesin e diagonalizimit të një matrice, gjetjen e bazës ortonormale të eigenhapësirave.  Ky \es një kurs i plotë i algjebrës lineare i cili vjen si një aplikim i një ushtrimi në dukje elementar atij të klasifikimit të prerjeve konike.  Në të vërtetë kjo \es edhe ana historike se si këto koncepte janë zhvilluar.  Kush i mëson këto koncepte pa kuptuar motivimin e tyre historik duhet të ndihet sikur ka parë një film që nuk e kuptoi dot kurrë.

Në Kreun~7  ne studiojmë ndërtimet gjeometrike.  Ndërtimet gjeometrike janë probleme klasike që mund të implementohen me sukses që në klasat e 5-ta dhe të 6-ta.  Nga ndërtimet gjeometrike mund të zhvillohet koncepti i fushës dhe shtrirjes algjebrike.  Ne japim një përshkrim të shkurtër të ndërtimeve algjebrike dhe vërtetimin e rezultateve kryesore. Ne supozojmë se lexuesi ka njohuri mbi fushat të krahasueshme me nivelin e \cite{alg}.

Ja vlen të përmendet se edhe një program ideal nuk do të funksiononte në se efektivi mësimdhënës nuk është në gjendje të kuptojë dhe implementojë me sukses këto koncepte intuitive. Pikërisht për këtë, formimi i një plejade matematikanësh me koncepte të qarta është më se i nevojshëm. 

Ne japim një sasi të konsiderueshme të literaturës, sidomos të disa botimeve të fundme në gjuhën Shqipe ku disa nga këto tema trajtohen në detaje \cites{alg, kalk, lin-alg}. 
Duke mos pretenduar se kjo \es fjala e fundit për mësimdhënien e matematikës, ne shpresojmë që mësuesit e matematikës do ti gjejnë këto ide një ndihmë në profesion e tyre  të vështirë. 

\section{Një vështrim i shkurtër mbi disa probleme të veçanta në programin e matematikës}

Njerëz të ndryshëm kanë perceptime të ndryshme mbi nivelin a nxënësve tanë sidomos para vitit 1990.  Pjesa më e madhe me të drejtë mendojnë se nxënësit tanë ishin të pregatitur mjaft mirë nga ana matematike.  Nga ana tjeter nuk mund të mos mohohet se programi i matematikës ka pasur dhe vazhdon të ketë probleme të shumta.  Ajo çka vihet re nga nxënësit shqiptarë që shkojnë me studime jashtë shtetit është se në shumicën e rasteve ata janë të pregatiur në mënyrë të pranueshme nga shkollat tona.  Ajo çfarë nuk kuptohet nga shumica e njerëzve është se nxënësit shqiptarë këtë e kanë arritur me një punë të jashtëzakonshme, se avantazhi më i madh i tyre karshi studentëve të tjerë ka qenë pikërisht disiplina në punë. Duke i shtuar kësaj edhe njohuritë e përfituara nga një sistem ku shkohej në shkollë gjashtë ditë në javë, ku bëhej matematikë çdo ditë për 12 vjet me rradhë, atëherë nxënësi i mirë i dalë nga shkollat tona kishte një avantazh të madh me shumicën e bashkëmoshatarëve te vet nga vendet e tjera.

Para shumë vjetësh një nxënës i vitit të tretë të një prej gjimnazeve të Vlorës, i bëri këtë pyetje njërit prej autorëve të këtij artikulli: 
 \textit{Profesor, më kanë thënë që shkolla në Amerikë është më e lehtë se në Shqipëri. Eshtë e vërtetë?}

Ka diçka me substancë në këtë pyetje.  Natyrisht ne po shohim përtej faktit që nxënës të shkëlqyer tanët shkojnë në shkolla mesatare në Amerikë apo vende të tjera. A janë programet e matematikës apo tekstet që përdoren më të thjeshta për nxënësin në vendet e tjera?  A mund të përmirësohen programet tona që mësimdhënia e matematikës të bëhet më efiçente?  Përgjigja \es "po".  Për nxënësin e mirë shkolla \es më e lehtë në Amerikë, sepse ka programe solide, tekste të shkëlqyera dhe pedagogë që kuptojnë jo vetëm subjektin që japin, por edhe shumë më tej.

Më poshtë po japim disa shembuj konkretë ku disa koncepte mund të trajtoheshin më mira në programin e matematikës. 

\begin{prob}
Gjeni bashkësinë e përcaktimit të funksionit
\[ f(x) = \frac 1 {\sqrt{x^2-4}} \]
\end{prob}

Shumica e nxënësve tanë të vitit IV, do ta bëjnë pa vështirësi këtë ushtrim, por problemi qëndron diku tjetër.  Çfarë është \textbf{bashkësia e përcaktimit}?  Ne përdorim termin \textbf{përkufizim} dhe jo \textbf{përcaktim}.  Ndryshimi midis dy termave në gjuhën Shqipe është i qartë.  Atëherë përse ky konfuzion në matematikë?  Ne i kemi ngjitur një koncepti fare të qartë matematik një term tepër konfuz.  Në fakt termi \textbf{bashkësia e përkufizimit} duhet përdorur sepse ne përkufizojmë një funksion nuk e përcaktojmë atë. 

Vazhdojmë me një shembull nga fizika.  

\begin{prob}
Një nga konceptet kryesorë në historinë e shkencës eshte koncepti i \textbf{shpejtësisë} (ose ajo çka në Anglisht quhet \textbf{velocity}).  Natyrisht koncepti i saktë i shpejtësisë jepet vetëm nëpërmjet derivatit.  
\end{prob}

A e kuptojnë nxënësit dhe mësuesit tanë ndryshimin midis \textbf{shpejtësisë} dhe \textbf{shpejtësisë mesatare} (në Anglisht \textbf{velocity} dhe \textbf{speed}). Në fjalët e një studente shqiptare që studion në Princeton, - "Ndryshimin e kuptova vetëm në Princeton".  

A kuptohet nga nxënësit dhe mësuesit tanë rëndësia historike e përkufizimit të saktë të shpejtësisë (pra të derivatit të funksionit).  Ishte pikërisht nevoja e përkufizimit saktë të shpejtësisë ajo që detyroi Newtonin së pari të zhvillonte Kalkulusin dhe së dyti të formulonte tre ligjet bazë te mekanikës të cilat ishin pikënisja e fizikës dhe e gjithë shkencës moderne. 

Shumica e gjuhëve sot kanë terma të veçantë për shpejtësinë dhe shpejtësinë mesatare, por në tekstet dhe programet tona këto ndërthuren si pa të keq.  Për shembull, kur një mësues i klasës së tretë jep problemin:\\

\textit{ Një udhëtar shkoi nga Vlora në Fier për 5 orë. Distanca Vlorë-Fier është 35 km.  Sa është shpejtësia?.} \\

Eshtë e qartë se këtu bëhet fjalë për shpejtësinë mesatare dhe jo shpejtësinë e çastit.  Madje në të folurën e përditshme shumica e përdorin fjalën \textit{shpejtësi} si shpejtësi mesatare. 
Janë pikërisht këto pakujdesi që gjenden kudo në tekstet tona që i bejnë gjërat tepër të vështira për nxënësit tanë. 

\subsection{Koncepti i dallorit të polinomeve}
Vazhdojmë me një shembull tjetër nga matematika e klasës së 10-të.  Dallori i polinomeve është një nga konceptet themelore të matematikës.  Në programet tona trajtohet fare cekët.  Supozojmë se 
u japin një grupi nxënësish dhe studentësh problemin e mëposhtëm:

\begin{pyetje}
Jepet polinomi quadratic 
\[ f(x) = ax^2 + bx + c\]
Përkufizoni dallorin e $f(x)$.  
\end{pyetje}

Shumica e nxënësve të shkollave të mesme si dhe gjithë mësuesët e matematikës do të japin përgjigjen e menjëhershme 
\[ \Delta = b^2 - 4 ac \]
A është kjo përgjigja e saktë? Shumë prej lexuesve do të thonë se 'po'. Por në se pyetja është të përkufizohet dallori i polinomit të gradës së tretë 
\[ f(x) = ax^3 + bx^2 + cx + d,\]
atëherë shumica dërrmuese e nxënësve dhe mësuesve nuk do të jenë në gjendje të japin një përgjigje të saktë. Kjo përforcohet edhe më shumë në se pyesim për dallorin e një polinomi të gradës $n\geq 4$. Pra ç'është dallori?  A e kuptojnë nxënësit dhe mësuesit tanë konceptin e dallorit?  

Në se në algjebrën e lartë jepet problemi që vijon:
\begin{prob}
Jepet një polinom $p(x)$ i gradës $n\geq 2$ 
\[ p(x) = a_n x^n + \dots a_1 x + a_0.  \]
Përkufizoni dallorin e k\"etij polinomi. 
\end{prob}
A është gati nxënësi të përgjithsojë konceptin e dallorit të polinomit kuadratik tek një polinom e gradës $n\geq 2$? 

Natyrisht, gjëja më e lehtë nga ana pedagogjike \es të flasësh për zgjidhjen e ekuacionit kuadratik kur $\D < 0$. Eshtë pikërisht këtu që në mënyrë fare të vazhdueshme mund të prezantohen numrat kompleksë dhe jo vetëm të prezantohen, por të shpenzohet kohë e mjaftueshme me bashkësinë e numrave kompleksë. 

\subsection{Një shembull nga trigonometria, rrethi njësi}

Nxënësit tanë kanë bërë trigonometri, në programin para 1990, për gjithë vitin III-të të shkollës së mesme.  Kjo siguronte një bazë solide për gjithë studentët e mirë. Në fakt në shumicën e programeve të matematikës në vendet perëndimore trigonometria bëhet vetëm një simestër ose ndërthuret në simestra të ndryshëm, por asnjëherë nuk zë kohën ekuivalente me dy simestra. 
Funksionet trigonometrike janë: 
\[ f(x) = \sin x,  \quad  f(x) = \cos x,  \quad f(x) = \tan x,  \quad f(x) = \cot x,  \]
Në rrethin njësi ato paraqiten si në figurë. 

\begin{minipage}{0.5\textwidth}
\centering
\begin{tikzpicture}[scale=1.9, cap=round]
  \def\costhirty{0.8660256}

  \colorlet{anglecolor}{green!50!black}
  \colorlet{sincolor}{red}
  \colorlet{tancolor}{orange!80!black}
  \colorlet{coscolor}{blue}

  \tikzstyle{axes}=[]
  \tikzstyle{important line}=[very thick]
  \tikzstyle{information text}=[rounded corners, fill=red!10, inner sep=1ex]

  \draw[style=help lines, step=0.5cm] (-1.2, -1.2) grid (1.2, 1.2);

  \draw (0, 0) circle (1cm);

  \begin{scope}[style=axes]
    \draw[->] (-1.3, 0) -- (1.3, 0) node[right] {$x$};
    \draw[->] (0, -1.3) -- (0, 1.3) node[above] {$y$};

    \foreach \x/\xtext in {-1,  -.5/-\frac{1}{2},  1}
      \draw[xshift=\x cm] (0pt, 1pt) -- (0pt, -1pt) node[below, fill=white]
      {$\xtext$};

    \foreach \y/\ytext in {-1,  -.5/-\frac{1}{2}, .5/\frac{1}{2},  1}
      \draw[yshift=\y cm] (1pt, 0pt) -- (-1pt, 0pt) node[left, fill=white]
      {$\ytext$};
  \end{scope}

\filldraw[fill=green!20, draw=anglecolor] (0, 0) -- (3mm, 0pt) arc(0:30:3mm);
  
\draw (15:2mm) node[anglecolor] {$\alpha$};

\draw[style=important line, sincolor] (30:1cm) -- node[left=1pt, fill=white] {$\sin \alpha$} +(0, -.5);

\draw[style=important line, coscolor] (0, 0) -- node[below=2pt, fill=white] {$\cos \alpha$} (\costhirty, 0);

\draw[style=important line, tancolor] (1, 0) --  node [right=1pt, fill=white]
    {
      $\displaystyle \tan \alpha \color{black}$
    } (intersection of 0, 0--30:1cm and  1, 0--1, 1) coordinate (t);

  \draw (0, 0) -- (t);

\end{tikzpicture}
\end{minipage}
\begin{minipage}{0.5\textwidth}

Një problem disi interesant nga klasa e 9-të është të shprehësh gjithë funksionet trigonometrike si funksione racionale të $\tan \frac {\alpha} 2$.  Një kërkesë e tillë i duket disi e çuditshme një nxënësi të vitit të III-të.  

\begin{prob}\label{trig}
Shprehni funksionet $\sin \a$ dhe $\cos \a$ në varësi të funksionit $\tan \frac {\alpha} 2$. 

\end{prob}

Shumica e nxënësve janë në gjendje ta bëjnë këtë ushtrim dhe të vërtetojnë formulat e mëposhtme:

\end{minipage}

\[  \sin \alpha = \frac    {2 \tan \frac {\alpha} 2} {1 + \tan^2 \frac {\alpha} 2}, \qquad \cos \alpha = \frac    {1 - \tan^2 \frac {\alpha} 2    } {1+ \tan^2 \frac {\alpha} 2} \]
Historia në programet tona mbyllet këtu, por kjo duhet të ishte vetëm fillimi i një dashurie të bukur midis nxënësit dhe trajetoreve algjebrike.  

Ne do të shohim më poshtë  se si një ushtrim kaq i thjeshtë mund të përdoret për të prezantuar nxënësit me ide dhe koncepte të thella matematike si ai i pikave racionale në trajektoret algjebrike.

Ndoshta një problem fare i natyrshëm që mund të bëhej me lehtësi në trigonometri ishte koncepti i rrenjëve të njësisë.  Për shembull, sa nga nxënësit apo mësuesët tanë mund ti përgjigjen saktë pyetje së mëposhtme.

\begin{pyetje}
Gjeni gjithë rrënjët e ekuacionit
\[ x^3 =1 \]
\end{pyetje}

Koncepti i rrënjëve të njësisë me shumë lehtësi na çon tek numrat kompleksë, shumëzimi i numrave kompleksë, grupi i rrethit njësi, dhe grupet ciklikë.

\subsection{Integrimi i funksioneve racionalë}

Le të supozojmë se audiencës sonë i shtrojmë pyetjen:
\begin{pyetje}
Jepet funksioni racional 
\[ f(x) = \frac {p(x) } {q(x)} ,  \]
ku $p(x), q(x)$ janë polinome me  koeficientë realë. Përshkruani një mënyrë për të llogaritur integralin
\[ \int f(x) \; dx \] 
\end{pyetje}
Ndoshta pyetja e mësipërme mund të thjeshtohet në një pyetje më konkrete,

\begin{exa} Njehsoni $\int \frac{2x^2-x+4}{x^3+4x}dx$.
\end{exa}

A është në gjendje maturanti ynë të mbaroj një ushtrim të tillë?  Për më tepër është një nga klasat më elementare të funksioneve.  Për një trajtim të kësaj teme lexuesi mund të lexojë \cite[Kap.~6]{kalk}.

Më poshtë po japin një zgjidhje të këtij ushtrimi për të kuptuar që trajtimi \es fare elementar.   Meqënëse $x^3+4x = x(x^2+4)$ dhe nuk mund të faktorizohet më tej mbi numrat realë, atëherë ne mund të shkruajmë
\[ \frac{2x^2-x+4}{x^3+4x} = \frac A x+ \frac{Bx+C}{x^2+4}\]
Duke shumëzuar të dy anët me $x(x^2+4)$,  kemi
\[ 2x^2-x+4 = A(x^2+4)+(Bx+C)x = (A+B)x^2+Cx+4A\]
Duke barazuar koefiçentët pranë fuqive të njëjta të $x$,  kemi $A=1$,  $B=1$,  dhe $C=-1$ dhe prej këtej
\[\int \frac{2x^2-x+4}{x^3+4x} \; dx=\int   \left( \frac 1 x+ \frac{x-1}{x^2+4} \right) \; dx\]
Në mënyrë që të integrojmë termin e dytë,  e ndajmë atë në dy pjesë:
\[\int \frac{x-1}{x^2+4} \; dx = \int \frac{x}{x^2+4} \; dx -\int \frac{1}{x^2+4}\]
Bëjmë zëvendësimin $u=x^2+4$ në integralin e parë të kësaj pjese,  kështu që $du = 2x\; dx$ dhe përfundimisht  kemi
\begin{equation*}
\begin{aligned}
 \int \frac{2x^2-x+4}{x^3+4x}dx & =\int \frac 1 x dx+\int \frac{x}{x^2+4}dx-\int \frac{1}{x^2+4}dx\\
& = \ln|x|+\frac 1 2 \ln(x^2+4)-\frac 1 2\tan^{-1}(x/2) + K\\
\end{aligned}
\end{equation*}

A mund të përgjithsohet kjo teknikë tek gjithë integralet e funksioneve racionalë?  Cila është teorema nga ana teorike që kjo metodë të funksionojë për gjithë funksionet racionalë?

Për çudi kjo metodë elementare nuk \es trajtuar në programet tona megjithëse bazohet në një fakt në dukje elementar, por shumë themelor nga ana teorike.

\begin{lem} 
Çdo polinom me koefiçenta realë faktorizohet si prodhim faktorësh linearë ose kuadratikë. 
\end{lem}

Arsyeja që kjo lemë \es e vërtetë \es sepse fusha e numrave realë \es invariant nën konjugimin kompleks. Pra një tjetër rezultat goxha i pafajshëm në dukje por me rendësi \te madhe teorike.

\subsection{Elipsi dhe rrethi}

Le të konsiderojmë disa ide nga gjeometria.  Jepet ekuacioni i nje elipsi
\[ \frac {x^2} {a^2} + \frac {y^2} {b^2} =1 \]

\noindent \begin{minipage}{0.45\textwidth}
\begin{tikzpicture}[scale=.75]
\draw[very thin,color=gray,step=1cm,dashed] (-3.5,-2.5) grid (3.5,2.5);
\draw[thick, ->] (-3.5,0) -- (3.5,0) node[below right] {$x$};
\draw[thick, ->] (0,-2.5) -- (0,2.5) node[left] {$y$};

\node[below right] at (3, 0) {$a$} ;
\node[above right] at (0, 2) {$b$} ;

\draw [-,samples=70,thick,red,domain=-3:3] plot(\x,  {(2/3)*sqrt(9-(\x)^2)} );
\draw [-,samples=70,thick,red,domain=-3:3] plot(\x,  {-(2/3)*sqrt(9-(\x)^2)} );
\end{tikzpicture}
\end{minipage}
\begin{minipage}{0.55\textwidth}

\begin{pyetje}
A mund të gjeni formula për sipërfaqen  dhe perimetrin e këtij elipsi?  
\end{pyetje}

Nxënësit e mirë padyshim që mund ta zgjidhnin një ushtrim të tillë.
Duke e zgjidhur këtë   ekuacion në lidhje me $y$,  ne marrim
\[ \frac{y^2}{b^2} = 1 - \frac{x^2}{a^2}=\frac{a^2-x^2}{a^2}\] 
ose
\[ y = \pm \frac b a \sqrt{a^2-x^2}. \]

Meqë elipsi është simetrik në lidhje me të dy boshtet koordinative,  sipërfaqja totale $A$ është katërfishi i sipërfaqes së kuadrantit të parë. 

\end{minipage}

Pjesa e elipsit në kuadrantin e parë jepet nga funksioni
\begin{equation}\label{funtion-el}
 y = \frac b a \sqrt{a^2-x^2} \ \ \ 0\leq x \leq a 
\end{equation}
dhe kështu
\[ \frac 1 4 A = \int_0^a \frac b a\sqrt{a^2-x^2} dx\]
Për të gjetur këtë   integral bëjmë zëvendësimin 
\[ \boxed{x=a\sin t}.\]
 Atëherë $dx = a \cos t \; dt$. Për të ndryshuar kufijtë e integrimit,  shohim se kur $x=0$,  $\sin t = 0$,  pra $t=0$ ndërsa  kur $x=a$,  $\sin t = 1$,  pra $t=\pi /2$. Gjithashtu
\[ \sqrt{a^2-x^2}=\sqrt{a^2-a^2\sin^2t}= \sqrt{a^2\cos^2 t} = a|\cos t|= a\cos t\]
meqë $0\leq t\leq \pi /2$. Prej nga,
\[
\begin{split}
A & = 4 \frac b a \; \int_0^a \sqrt{a^2-x^2}  \; dx  = 4 \frac b a \; \int_0^{\pi /2}a\cos t \cdot a\cos t \; dt  \\
& = 4 a b \int_0^{\pi /2}\cos^2t \; dt  = 4 ab\int_0^{\pi /2}\frac 1 2 \left(1+\cos 2t \right) \; dt\\
& =\left. 2ab \left[  t+\frac 1 2 \sin 2t  \right]    \right|_0^{\pi/2}  = 2ab \left(\frac{\pi}{2}+0-0 \right) = \pi ab.\\
\end{split}
\]
Vërtetuam  kështu se sipërfaqja e elipsit me gjysmëboshte $a$ dhe $b$ është 
\[ \boxed{ A= \pi ab }\]
 Në veçanti kur marrim $a=b=r$,  kemi vërtetuar formulën se sipërfaqja e rrethit me rreze $r$ është $\pi r^2$.

Në nivelin e shkollave tona ku lloj integrali është e mundur që të bëhet. Ajo çka është e vecanta është se na orienton në dy probleme interesante.  Së pari, përgjithësimi i zëvendësimit të mësipërm na çon tek \textit{zëvendësimet trigonometrike} dhe prej andej tek formulat e integrimit për 
\[ \int \sqrt{a^2-x^2 } \; dx, \quad \int \sqrt{a^2+x^2 } \; dx, \quad \int \sqrt{x^2-a^2 } \; dx\]
Së dyti, gjetja e perimetrit të elipsit na çon drejt një lëndine plot me lule, atë të integraleve eliptike. 

Gjatësia e grafikut të një funksioni $y=f(x)$, \; $a \leq x \leq b$ jepet nga formula
\[ L = \int_a^b \sqrt{1+ \left[ f^\prime (x) \right]^2 } \; dx, \]
shihni \cite[Kapitulli 9]{kalk}. Në rastin tonë, $f(x)$ jepet ne Eq.~\eqref{funtion-el}. Derivati i tij është 
\[ f^\prime (x) = \frac b {2a } \;   \frac 1 {\sqrt{a^2 - x^2}} \; (-2x) = - \frac b a \, \frac x {\sqrt{a^2 - x^2}} \]
Pra perimetri i elipsit është
\[ P= 4 \int_0^a \sqrt{ 1  + \frac {b^2} {a^2}   \cdot  \frac {x^2} {a^2 - x^2} } \; dx = 4 \; \int_0^a \sqrt{ \frac {a^2-e^2x^2} {a^2-x^2}   } \]
ku $e= \frac {a^2-b^2} {a^2}$. Duke bërë zëvendësimin 
\[ \boxed{ x=a\sin t} \]
lexuesi të vërtetojë se 
\[ P= 4 a \int_0^{\frac {\pi} 2} \sqrt{1 - e^2 \sin^2 t } \; dt \]
Ky quhet një \textbf{integrali i plotë eliptik i llojit të dytë}.  Prej tij marrin emrin \textit{trajektoret eliptike} që janë nga objektet mjaft popullore të matematikës.  Integralet eliptike ishin edhe pikënisja e punës së Abelit dhe Jacobit që vazhdoi më pas me Riemann.

Eshtë për të ardhur keq që në programet e shkollave të mesme integralet eliptike nuk trajtohen fare. Ato janë të vështira për tu zgjidhur, por një pjesë e konsiderueshme e matematikës moderne nisi pikërisht nga integrale të tilla si për shembull funksionet theta \cite{MR3525572}, integralet hypereliptike dhe supereliptike, trajektoret hipereliptike dhe supereliptike \cite{MR3525570}, etj.

Më poshtë po japim një tjetër problem në dukje shumë afër problemit të mësipërm, por që na çon në një fushë tjetër të matematikës.

\begin{pyetje}
A mundet qe ky elips të transformohet në një rreth me një ndryshim koordinatash të planit?  A ndryshojnë sipërfaqja dhe perimetri gjatë këtij transformimi?
\end{pyetje}

Zgjidhja e këtij problemi \es një mundësi për të prezantuar nxënësit me transformimet lineare dhe matricat ortogonale. Ne flasim pak për to në vijim. 

\subsection{Prerjet konike}

Kjo ideja e transformimeve na çon edhe më thellë. Naturisht në se ne transformojmë planin në një sistem tjetër koordinativ, për shembull
origjina shkon tek pika $(2, 3)$, atëherë ekuacioni i elipsit bëhet 
\[ \frac {(x-2)^2} {a^2} + \frac {(y-3)^2 } {b^2} =1. \]
Ekuacion që ndryshe shkruhet 
\[  b^2 x^2 + a^2 y^2 - 2b^2 x - 6 a^2 y = a^2 b^2 - 4 b^2 - 9 a^2 \]
Ne e dimë qe ky është ekuacioni i një elipsi, edhe pse kjo nuk është plotësisht e lehtë të verifikohet.  

A mund ta përgjithsojmë këtë rast në një teori më të përgjithshme.  Për shembull, jepet ekuacioni i përgjithshëm i gradës së dytë
\[ ax^2 + b xy + cy^2 + dx + fy + s=0  \]
\begin{pyetje}
Për ç'vlera të koefiçentëve $a, b, c, d, f, s$ grafiku i mësipërm është një rreth, elipse? 
\end{pyetje}

Por ndoshta edhe më themelore \es pyetja

\begin{pyetje}\label{pyetje:conic}  Jepet ekuacioni i gradës së dytë 
\begin{equation}\label{conic}
  ax^2 + b xy + c y^2  +    d x + e y=    \l 
\end{equation}  
Çfarë paraqet një ekuacion i tillë nga ana grafike?  A mund të gjendet një metodë që kësaj pyetje t'i jepet përgjigje thjesht nga koefiçentët $a, b, c, d, e, \l$, pa ndërtuar grafikun? 
\end{pyetje}

Nxënësi ynë e di se grafikët e ekuacineve të tilla janë grafikët e prerjeve konike, pra prerjet e një koni të dyfishtë me një plan si në figurën më poshtë. 
\begin{figure}[hbp]
\begin{center}
\includegraphics[scale=0.33]{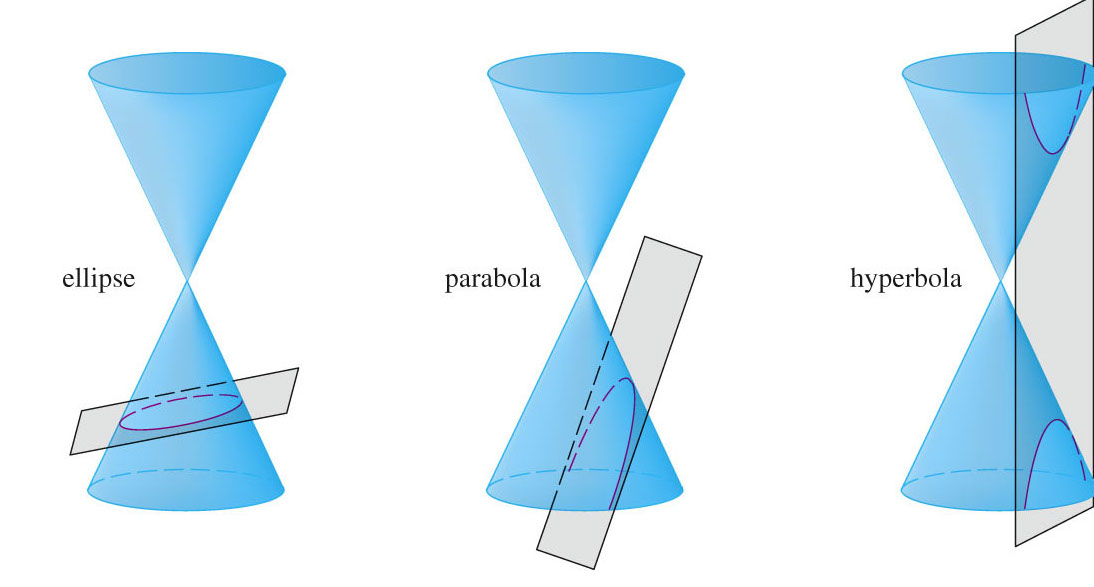}
\caption{Prerjet konike}
\label{10p670a}
\end{center}
\end{figure}

Nxënësit e shkollave tona e dinë që, elipsi, parabola, dhe hiperbola kanë përkatësisht ekuacione
\begin{equation}\label{conics}
 \frac {x^2} {\l_1^2} + \frac {y^2} {\l_2^2} =1, \qquad y^2=4\l_1 x, \qquad \frac {x^2} {\l_1^2} - \frac {y^2} {\l_2^2} =1
\end{equation} 
për disa vlera jozero $\l_1, l_2$. Rrethi \es rasti i veçantë i elipsit kur $\l_1=\l_2$.  

Secili prej këtyre ekuacioneve \es shumë herë më i thjeshtë se ekuacioni Eq.~\eqref{conic} dhe kjo sepse sistemi koordinativ \es zgjedhur në qendrën e rrethit, elipsit, hiperbolës, apo kulmin e parabolës dhe boshtet në mënyrë të përshtatshme.  Këto ekuacione kanë të mirën e madhe që ne e dimë formën e grafikut thjesht duke parë ekuacionin.  Atëherë shtrohet pyetja, e kuptueshme nga çdo gjimnazist:

\begin{pyetje}
Si duhet ndryshuar sistemi yne koordinativ që ekuacioni \eqref{conic} të transformohet në një nga ekuacionet e mësipërme?
\end{pyetje}

Rëndësia e kësaj pyetje dhe konceptet që mund të zhvillohen duke ju përgjigjur kësaj pyetje janë një mrekulli e vërtetë.  Për herë të parë nxënësi ekspozohet tek ideja që ekuacionet janë diçka pak e rëndësishme, diçka që ndryshon.  Eshtë objekti gjeometrik ai që nuk ndryshon.  Apo ka vallë ndonjë koncept sasior që mund të shprehet në varësi të koefiçentëve dhe që nuk ndryshon pavarësisht ndryshimit të sistemit koordinativ?  

Si pa kuptuar ne jemi futur në teori invariantesh, një nga degët më aktive të matematikës së shekullit XIX. Gjithashtu, pikërisht ky diskutim kaq modest është fillimi i gjeometrisë ose më mirë i gjeometrisë algjebrike, një nga degët më aktive të matematikës moderne.  Ne do ti japim zgjidhje pyetjes së mësipërme ne kreun në vazhdim.

\begin{pyetje}
A mund të gjejmë një parametrizim racional të prerjeve konike?  Ose me konkretisht, a mund të shprehet ekuacioni i çdo prerje konike si ekuacion parametrik $(x(t), y(t) )$ ku $x(t)$ dhe $y(t) $ janë funksione racionalë? 
\end{pyetje}

Kujtoni që \es një problem elementar që zakonisht \es bërë  me nxënësit e mirë që ekuacionet parametrik i   elipsit \es 
\[ x= a\cos \theta, \qquad   y = b \sin \theta \]
i parabolës
\[ x= at^2,   \qquad  y= 2at \]
dhe i   hiperbolës
\[ x= a \sec \theta, \qquad y = b \tan \theta . \]
Natyrisht parabola ka një parametrizim racional.  Për elipsin dhe hiperbolën mjafton të zevendësojmë $t= \tan \frac {\theta } 2$ dhe ne bazë të Problemit~\ref{trig} ne marrim një parametrizim racional.  Pra kemi se prerjet konike gjithmonë kanë një parametrizim racional.  Për entuziastuet e gjeometrisë algjebrike, kjo do të thotë se 
\textbf{prerjet konike janë trajektore algjebrike me genus zero}.  
Kjo konkluzion   na hap një pyetje tjetër:

\begin{pyetje}
A ka ndonjë rëndësi ky parametrizimi racional?  
\end{pyetje}

Si pa dashur jemi futur në problemin klasik të Diofantit, atë të gjetjes së zgjidhjeve integrale (në mënyrë ekuivalente racionale për ekuacionet homogjenë) të ekuacioneve algjebrikë.  

Ne nuk do të zgjerohemi shumë në këtë fushë, por thjesht duam të theksojmë se gjithë teoria moderne e numrave vërtitet rreth këtij problemi.  Duam gjithashtu të theksojmë që një ekuacion algjebrik që ka një parametrizim racional ka një bashkësi të pafundme zgjidhjesh racionale. 

Eshtë gjithashtu me vlerë për tu theksuar se shumica e ekuacioneve polinomiale kanë vetëm njue numër të fundëm zgjidhjesh racionale. Kjo \es teorema Faltings (83) dhe konsiderohet si teorema e shekullit të XX. 

Pra prerjet konike nga ky këndvështrim duken mjaft speciale dhe nga një ushtrim elementar i trigonomestrisë ne kuptojmë përse ato janë tue tilla. 

\subsection{Maximumet dhe minimumet lokale te funksioneve}
Vazhdojmë me një shembull tjetër nga analiza.  Një pjesë e mirë e vitit IV në analizë kalohet me gjetjen e maksimumeve dhe minimumeve lokale të funksioneve me një ndryshore.  Nxënësit harxhojnë me javë të tëra duke bërë probleme optimizimi.  Asnjë fjalë në analizën e vitit IV nuk thuhet për funksionet me dy ndryshore.  Mendoni sa kuriozitet do të zgjonte tek nxënësi ideja e gjetjes se majave të kodrave dhe thellësive të detit, apo llogaritja e sipërfaqes së një relievi.  Ishtje një pyetje që shqetësonte autorin e dytë të këtij artikulli si adoleshent: \textit{Sa sipërfaqe reale tokë ka Shqipëria? }

Pra problemi i mëposhtëm mund ë shtrohet tek maturantët tanë dhe ndoshta të punohet vetëm me ata më të talentuarit.
\begin{prob}
Jepet një sipërfaqe 
\[ z = f(x, y) \]
e përkufizuar në një zonë $D \subset \R^2$. Të gjenden maksimumet dhe minimumet lokale të kësaj sipërfaqeje. 
\end{prob}

Natyrisht analiza \es një minierë e vërtetë idesh dhe aplikimesh.  Gjetja e vorbullës së një lumi, volumi i një kodre, kurbatura e një kthese, shpejtësia dhe nxitimi i një objekti që lëviz në hapësirë, ligjet e Keplerit, etj.  Duhet pranuar se analiza e funksioneve me shumë ndryshore edhe në programet tona universatere \es bërë cekët dhe ilustrimi me aplikime nga fizika dhe inxhinjeritë gjithmonë ka qenë i varfër. 

Gjithashtu ja vlen të përmendet se shumë prej këtyre temave trajtohen në shkollat e mesme në shumë vende. Për shembull në Amerikë nxënësit që janë të orientuar drejt shkencave dhe inxhinjerive futen ne klasa të avancuara dhe marrin lendë si Kalkulusi I dhe Kalkulusi II që në shkollë të mesme.

\subsection{Ndërtimet gjeometrike}

Ata që kuptojnë historinë e matematikës të paktën njëherë në jetën e tyre duhet të shkojnë të vizitojnë  varrezat Albani \ne Gottingen ku pushon një prej njerëzve më me influencë në historinë e njerëzimit, Carl Friedrich Gauss.  Gauss \es i njohur për shumë rezultate të famshme \ne matematikë dhe një prej tyre ishte edhe ndërtimi me vizore dhe kompas i një 17-këndëshi \te rregullt ose siç quhet ndryshe \textit{heptadecagon}, pas 2000 vjetësh përpjekje nga matematikanë të shumtë. Ishte një nga arritjet që e bënin Gaussin krenar më shumë se çdo gjë tjetër.  
Çfarë \es kaq e vështirë për këto ndërtimet gjeometrike që i dha Gausit kaq famë dhe krenari?

\begin{figure}[htbp] 
   \centering
   \includegraphics[width=4in]{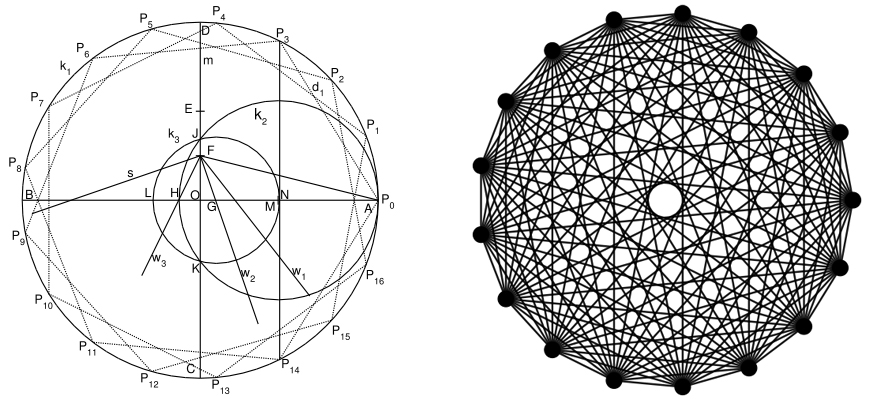} 
   \caption{Heptadecagon ose 17-këndëshi i rregullt }
   \label{fig:example}
\end{figure}

Në tekstet tona të shkollës 8-vjeçare dhe të mesme, gjithmonë ka pasur copëza historike mbi ndërtimet gjeometrike, por këto nuk janë shfrytëzuar për të ndërtuar koncepte matematike mbi to.  Një libër i shkëlqyer i Petro Priftit, \textit{Probleme të zgjidhura të gjeometrisë} ka pasur plot probleme të vlefshme, por ky tekst ishte gjithmonë diçka që shfrytëzohej vetëm nga studentët e talentuar. Pak përdoreshin probleme të tilla për të orientuar nxënësin drejt koncepteve të matematikës bashkëkohore, madje shumica e mësuesve nuk ishin fare në dijeni të këtij libri kaq të vyer. 
Le të përpiqemi ti organizojmë disi këto probleme klasike të gjeometrisë.

Në Greqinë e lashtë kishte   disa problema klasike. Këto probleme janë nga  gjeometria në natyrë dhe përfshijnë ndërtimet vetëm  me  vizore  dhe kompas.   Problemet mund të formulohen si më poshtë.

\begin{enumerate}
\item Jepet një kënd i \c cfarëdoshëm,  a mund ta ndajmë këndin në tre kënde të barabartë duke përdorur vetëm  vizore dhe kompas?

\item Jepet një rreth i \c cfarëdoshëm,  a mund të ndërtojmë një katror me të njëjtën sipërfaqe,  duke përdorur vetëm  vizore dhe kompasin?

\item Jepet një kub,  a mund të ndërtojmë brinjën e një kubi tjetër,  i cili të ketë dyfishin e vëllimit të kubit origjinal duke përdorur  vetëm  vizore  dhe kompas.
 
\item  Për cilat  $n$,   $n$-këndëshi i rregullt është  i ndërtueshëm?
\end{enumerate}

Pas   përpjekjeve  dymijë vje\c care nga ana e matematikanëve,  u tregua,    ndërtimet në tre problemet e para janë të    pamundura.     Në Kreun~6 ne do të japim disa ide se si këto tema mund të trajtohen që në shkollë të mesme.


\section{Historia e shkencës - një mjet kritik dhe konstruktiv për programin e matematikës}

Ka mendime dhe këndvështrime të ndryshme për pregatitjen e nxënësve në matematikë.  Këto pikëpamje kanë lidhje me shoqërinë, traditën e një vendi, filozofinë e shkollës përkatëse.  Ajo që pothuajse \es e pranueshme nga të gjithe \es \qe intuita matematike zhvillohet në një moshë të vogël dhe konceptet merren herët.  Diskutimet sot në ambjentet pedagogjike janë kryesisht mbi rrugët që duhen ndjekur me moshat e reja.  

Shkolla jonë ka të theksuar përsëritjen e vazhdueshme dhe për një kohë të gjatë dhe kjo shpesh \es bërë në kurriz të anës intuitive.  Po të kemi parasysh se para viteve 1990 fëmijët shkonin 6 ditë \ne shkollë dhe bënin përditë matematikë, atëherë duhet të pranojmë se femijët shqiptarë harxhonin një pjesë të konsiderueshme të fëmijërisë me matematikën.  Po të shohësh materialin që bëhej për 12 vjet shkollë \es më pak se shumica e vendeve të tjera.  Ky material përsëritej aq herë sa bëhej i mërzitshëm, sidomos për nxënësit e mirë.  
Ç'mund të bëhet me këta nxënës?  Si mund të zhvillohet intuita tek të gjithë nxënësit dhe veçanërisht tek nxënësit me prirje në matematikë?  

Zakonisht rruga më e mirë \es kur konceptet zhvillohen në mënyrën në të cilën janë zbuluar.  Pra shtrohet pyetja, orientohen nxënësit drejt ideve të dukshme edhe pse ato mund të jenë të gabuara. Eshtë e nevojshme \qe nxënësi të kuptojë se çfarë nuk funksionon dhe të detyrohet vetë, ndoshta edhe duke vuajtur pak, që të gjejë zgjidhjen e duhur.  Faktet historike, vënia e çdo teoreme në ambjentin historik \es një ndihmë shumë e madhe.  Historia e matematikës moderne \es një dramë e ngjeshur emocionale e shekujve XVIII, XIX, XX.  

Edukimi i përgjithshëm i matematikës sot \es distancuar nga kjo lloj mësimdhënie.  Ne shkojmë në klasë me një tufë teoremash dhe faktesh dhe detyrojmë nxënësit ti mësojnë ato.  Një mënyrë tepër komode për mësuesin mediokër, por tepër katastrofike për nxënësin.  Megjithatë ka plot raste dhe shkolla ku matematika mësohet ndryshe.  Ne do ti sygjeronin lexuesit librin \textit{Perfect rigor} \cite{MR2598223}, mbi jetën e një prej matematikanëve më të medhenj të këtij shekulli Grigory Perelman.  Nje pjesë e mirë e librit flet për \textit{Një shkollë të mrekullueshme} ku Grigory 12-vjeçar trajnohej përditë.  

Një shkollë tjetër e nisur vitet e fundit \es \textit{The Proof school} \ne San Francisco, ku intuita dhe konceptet klasike të matematikës marrin përparësi.

Natyrisht vështirësia më e madhe në këtë lloj shkolle \es mungesa e mësuesve të aftë, të cilët të kenë njohuri të thella të matematikës së lartë dhe të dinë se si këto koncepte ti ingranojnë në programet mësimore. 
Matematika, kur mësohet nga njerëzit e duhur hap shumë dyer dhe labirinte të thella në mendjen e një fëmije.  Problemi \es se si mund ti gjejmë këta vizionerë, këta njerëz të duhur që kuptojnë se ç'donte të bënte Gauss, Abel, Jakobi, që njohin thellë matematikën e shekullit XIX, që janë gati të sakrifikojnë karierën e tyre shkencore për tu mësuar disa 12-vjeçarëve matematikë.  Këta nuk \es e lehtë ti gjesh! Prandaj shoqëria duhet të bëj një përpjekje të madhe që të rekrutojë dhe pregatisë njerëz të tillë.

Për mësuesit e apasionuar të cilët duan të gjejnë rrugë të reja se si historia e matematikës mund të futet në klasën e matematikës ne do t'u sygjeronim librat e Felix Klein 
\cites{MR0076336, MR0015349, MR0055397, MR0080930,  MR549187, MR2098410, MR3495524,  MR3495525}, Wilder \cite{MR0297501}, Roberts \cite{MR767867}, dhe veçanërisht librin nga Arnold \cite{MR3409220}, 

\section{Ekuacionet algjebrike}

Tani japin disa ide mbi disa tema që janë themelore në historinë e zhvillimit të matematikës dhe që mund të futen me sukses ne programet e shkollave të mesme.

\subsection{Ekuacionet me grade 2 dhe 3}

Zgjidhja e ekuacioneve polinomiale ka qenë dhe vazhdon të jetë një nga problemet themelore të matematikës që ka nxitur zhvillimin e disa prej degëve më elegante dhe më produktive të matematikës si teoria Galua, gjeometria algjebrike, gjeometria Diofantine, etj. Më poshtë japim një ide se si zgjidhen disa nga ekuacionet me gradë të vogël me një ndryshore. 
Të gjithë polinomet kanë koefi\c cientë në  $\C$.  

Kur kemi të bëjmë me një ekuacion me një ndryshore
\begin{equation}\label{eq-1}
x^n\ +\ a_{n-1} x^{n-1}\ +\ \dots \ + \ a_1 x\ + \ a_0 \ = \ 0
\end{equation}
thjeshtimi i parë  është të zëvendësojmë
\[ x\ = \ y\ -\ \frac{a_{n-1}}{n}\]
i cili rezulton në një ekuacion
\[ y^n\ +\ b_{n-2} y^{n-2}\ +\ \dots \ + \ b_0\ = \ 0\]
me termin $y^{n-1}$ zero. Ky e zgjidh ekuacionin kuadratik (duke plotësuar katrorin). Pra  për $n=2$,   ne thjesht marrim $y^2=-b_0$.

Konsiderojmë tani rastin për $n=3$,  pra,  ekuacionin
\begin{equation}\label{eq-2}
y^3\ +\ ay\ + \ b\ = \ 0.
\end{equation}
Duke zëvendësuar $y=u+v$,   marrim
\begin{equation}\label{eq-3}
u^3\ +\ v^3\ +\ 3\ \left( uv+\frac{a}{3} \right)\ (u+v) \ + \ b\ = \ 0
\end{equation}
qe është i vërtetë,  në qoftë se $u$ dhe $v$ kënaqin
\begin{equation}
u^3\ +\ v^3\ \ = \ -b,  \ \ \ \ \ uv\ \ = \ -\frac{a}{3}
\end{equation}
Ekuacioni i fundit na jep
\[ (Z-u^3) \ (Z-v^3)\ = \ Z^2\ +\ bZ\ -\ \frac{a^3}{27}\]
pra,  
\begin{equation}
\begin{split}
u^3\  & = \ \frac{-b\ +\ \sqrt{b^2\ +\ \frac{4a^3}{27}}}{2}\\
v^3\ & = \ \frac{-b\ -\ \sqrt{b^2\ +\ \frac{4a^3}{27}}}{2} \\
\end{split}
\end{equation}

Nga ana tjetër,  ne mund të zgjedhim $u$ dhe $v$ si rrënjë kubike të përshtatshme të anës së djathtë,  të
tilla që,  Eq.~\eqref{eq-3}  të jetë i vërtetë. Atëherë $u+v$ është një nga zgjidhjet e Eq.~\eqref{eq-2},  të tjerat i marrim
nga zgjedhjet e ndryshme të rrënjëve kubike. Më saktë,  në qoftë se,  $\epsilon$ është një rrënjë primitive e tretë e njëshit, pra $\epsilon^3 =1$,   atëherë
zgjidhjet e Eq.~\eqref{eq-1}, janë 
\[y_1= \ u\ +\ v,  \ \ \ y_2=\ \epsilon u\ +\ \epsilon^2 v,  \ \ \ Y_y=\ \epsilon^2 u\ +\ \epsilon v, \]
të cilat mund të verifikohen duke faktorizuar  $(y-y_1)(y-y_2)(y-y_3)$. Këto janë \textbf{formulat Kardano}, 
të cilat zakonisht shkruhen si
\begin{equation}
 y_i\ = \ \left(-\frac{b}{2}\ +\ \sqrt{\frac{b^2}{4}\ +\ \frac{a^3}{27}} \right)^{1/3} \ +\ \left(-\frac{b}{2}\ -\ \sqrt{\frac{b^2}{4}\ +\ \frac{a^3}{27}} \right)^{1/3}
\end{equation}

Shohim se,   këto formula kanë disa simetri,  duke filluar nga zgjedhjet e ndryshme të rrënjëve katrore,  të rrënjëve kubike dhe  gjithashtu nga zgjedhja e një rrënje të tretë të njëshit.
Gjëja pozitive  është,   që \c cështja e simetrive midis zgjidhjeve mund të përkufizohet pa përdorur formulën e shtjellur të zgjidhjeve. 
 
Ideja themelore e teorisë Galois tani mund të formulohet thjesht  si vijon: \\
\textit{Për çdo $n\geq 4$,   zëvendësojmë formulën e shtjellur për  zgjidhjet (megjithëse  ajo nuk ekziston) duke përdorur grupin Galois}.

Për më tepër për zgjidhjen e ekuacioneve dhe një prezantim mbi teorinë e Galois shihni \cite{alg}. 

Në fakt koncepti i grupit në algjebrën moderne lindi pikërisht nga koncepti i bashkësisë së simetrive midis rrënjëve të një polinomi.  Megjithë rolin qëndror të teorisë Galua në algjebër, në programet tona, duke përfshirë edhe ato universitare, kjo teori as përmendej fare.  Pra studentët tanë të shkencave asnjëherë nuk arritën të bëjnë lidhjen apo të kuptojnë motivimin se përse studioheshin grupet apo fushat në universitet.  Dhe kur ata vetë nuk e kuptonin këtë lidhje nuk mund të pretendosh që ata t'u shpjegonin nxënësve në klasën e X degëzimet apo idetë që dalin nga zgjidhja e ekuacionit.

\section{Diskriminanti dhe identitetet e Newtonit}
Në proçesin e diskutimit të rrënjëve të një polinomi është shumë i rëndësishëm fakti në se ndonjë prej këtyre rrënjëve përsëritet.  Në fund të fundit, ekuacioni
\[ (x-2)^4 \, (x-3)^5=0 \]
është më i lehtë për tu zgjidhur edhe pse është i gradës 9.  Pyetje e mëposhtme është fare e natyrshme për çdo nxënës së klasës IX apo X.

\begin{pyetje}
Jepet polinomi i gradës $n\geq 2$
\[ 
f(x)\ \ = \ a_n  x^n\ +\ a_{n-1}x^{n-1}\ +\dots +\ a_0 \ \ = \ a_n (x-\a_1) \dots (x-\a_n)
\]
ku $a_n \neq 0$. 
A mund të gjeni një kriter për koefiçentët $a_0, \dots , a_n$ që $f(x)$ të ketë rrënjë që përsëriten?
\end{pyetje}

Ky nuk është një problem i vështirë, por na çon në një nga konceptet më të rëndësishëm të matematikës, atë të \textbf{dallorit} ose \textbf{diskriminantit}. Le ti rradhisim gjithë rrënjët si më poshtë
\[ \a_1, \a_2, \dots , a_{n-1}, \a_n \]
dhe marrim prodhimin 
\[
\D_f \ = \  \Pi_{i \neq j} (\a_i - \a_j) 
\]
Natyrisht, kemi një factor $(-1)$ që varet nga renditja e rrënjëve. Për ta bërë $\D_f$ të pandjeshëm nga kjo renditje e rrënjëve si edhe nga koefiçenti udhëheqës $a_n\neq 0$ ne modifikojmë përkufizimin si më poshtë: 
\begin{equation}\label{disc}
\D_f \ = (-1)^{n(n-1)/2} \, \cdot \, a_n^{2n-2} \, \cdot \,  \Pi_{i \neq j} (\a_i - \a_j) 
\end{equation}
Atëherë, lema e mëposhtme është e besueshme për çdo nxënës 
\begin{lem}
$f(x)$ ka rrënjë që përsëriten \att $\D_f = 0$. 
\end{lem}
Pikërisht, është $\D_f$ në Eq.~\eqref{disc} se si duhet të përkufizohet discriminanti i çdo polinomi.  Ky është koncepti i natyrshëm, koncepti që mbahet mend, dhe koncepti që na tregon vetitë e diskriminantit. 

Megjithatë ne ende nuk i jemi përgjigjur pyetjes së mësipërme.  Ne nuk i njohim rrënjët e $f(x)$.  Pra, a mund të themi diçka për shumëfishmërinë e rrënjëve pa i gjetur rrënjët? Kjo pyetje është ekuivalente me

\begin{pyetje}
Jepet polinomi 
\[ 
f(x)\ \ = \ a_n  x^n\ +\ a_{n-1}x^{n-1}\ +\dots +\ a_0 \ \ = \ a_n (x-\a_1) \dots (x-\a_n)
\]
A mund të shprehet diskriminanti $\D_f$ në varësi të  koefiçentëve $a_0, \dots , a_n$?
\end{pyetje}

Ushtrimi i mëposhtëm është një ushtrim që çdo nxënës duhet ta ketë bërë të paktën një herë në jetë.

\begin{exe}
Jepet polinomi quadratik
\[f(a)\ = a \ x^2\ +\ bx\ +\ c \ = \ a \ (x-\a_1)\ (x-\a_2).\] 
Atëherë, nga përkufizimi dallori është 
\[ \D_f =  a^2 \, (\a_1 - \a_2)^2 \]
Vërtetoni se  
\[ \D_f\ = \ b^2\ - \ 4 a c.\]
\end{exe}

Natyrisht, zgjidhja e ushtrimit të mëposhtëm kërkon vetëm njohuri elementare të algjebrës dhe mund të bëhet nga shumica e nxënësve të klasës X. 

\begin{exe}
Jepet polinomi kubik
\[f(a)\ = a x^3 + b x^2  + c x + d \] 
Atëherë, vërtetoni se  
\[ \D_f\ = b^2c^2-4ac^3-4b^3d-27a^2d^2+18abcd.\,\]
\end{exe}


Tani le të përpiqemi të nxjerrim  një formulë të shtjellur  për discriminantin.   Le të jepet  $f(x)$ si vijon
\[ 
f(x)\ \ = \ \ x^n\ +\ a_{n-1}x^{n-1}\ +\dots +\ a_0 \ \ = \ \ (x-\a_1) \dots (x-\a_n)
\]
dhe përkufizojmë
\[ \D \ \ = \ \ \prod_{i \neq j} \ (\a_i- \a_j)
\]
Matrica
\[
X := \ \ \
\begin{pmatrix}
1 & \dots  & 1\cr \a_1& \dots  & \a_n\cr. & \dots  &. \cr. & \dots  &. \cr. & \dots  &. \cr \a_1^{n-1}& \dots  & \a_n^{n-1}
\end{pmatrix}
\]
ka
\[ \det(X) \ \ \ = \ \ \ \prod_{i>j} \ (\a_i- \a_j)
\]
nga formula e mirënjohur e përcaktorit të matricës Vandermonde.

\begin{exa}  Vërtetoni formulën
\[ \D  =   \det(X)^2 =  \det(X X^t) =  \det
\begin{pmatrix}
S_0 & S_1 & \dots  & S_{n-1}\cr S_1& S_2& \dots  & S_{n}\cr
. &. & \dots  &. \cr
. &. & \dots  &. \cr. &. & \dots  &. \cr
S_{n-1}& S_{n} & \dots  & S_{2n-2}
\end{pmatrix}
\]
ku 
\[ S_\mu\ := \ \ \ x_1^\mu\ +\ \dots \ +\ x_n^\mu.\]

\proof 
Duhet të shprehim shumat fuqi $S_\mu$ në lidhje me  funksionet simetrikë elementarë
$$\s_\nu\ \ = \ \ \sum_{i_1<i_2< \dots <i_\nu}\ x_{i_1}x_{i_2}\ \dots \ x_{i_\nu}$$
(ku vetia bazë e tyre është,   që $\s_\nu(x_1,  \dots , x_n)\ = \ (-1)^{\nu}\ a_{n-\nu}$). Kjo mbështetet te
 \textbf{ identitetet e Newtonit} (c.f. Cox et al.,  p. 317):
\[
\begin{split}
& S_\mu  - \s_1 S_{\mu-1} + \dots  + (-1)^{\mu-1}\s_{\mu-1}S_1 + (-1)^{\mu}\mu \s_{\mu}= 0, \ \ \mbox{for } 1\le\mu\le n\\
& S_\mu - \s_1 S_{\mu-1} + \dots  + (-1)^{n-1}\s_{n-1}S_{\mu-n+1}  + (-1)^{n}\s_{n} S_{\mu-n} = 0, \ \ \ \mbox{for}\ \mu> n
\end{split}
\]
\end{exa}

Jepet $z$ një variabël i ri,   përkufizojmë
\[ \s(z)\ \ = \ \ \prod_{i=1}^n\ (1\ - \ x_iz)\]
Atëherë 
\begin{equation}
\begin{split}
& \frac{-z\s'(z)}{\s(z)}\ \ = \ \ \frac{z\ \sum_{i=1}^n\ x_i \prod_{jnë i}\ (1\ - \ x_jz)}{\s(z)}\
\ =
 \ \ \sum_{i=1}^n\ \frac{x_iz}{1\ - \ x_iz} \\
& = \,  \sum_{i=1}^n\ \sum_{\nu=1}^\infty\ x_i^\nu z^\nu\ \ =
 \ \ \sum_{\nu=1}^\infty\ (\sum_{i=1}^n\ x_i^\nu)\ z^\nu\ \ =
 \ \ \sum_{\nu=1}^\infty\ S_\nu\ z^\nu
\end{split}
\end{equation}
Kështu që,  ne marrim identitetin në vazhdim ndërmjet serive fuqi formale në $z$:
$$ \s(z)\ \sum_{\nu=1}^\infty\ S_\nu\ z^\nu \ \ = \ \ -z\s'(z)$$
Vetia bazë e funksioneve simetrikë elementarë na \c con në
$$\s(z)\ \ = \ \ \sum_{\mu=0}^n\ (-1)^\mu\s_\mu z^\mu$$
Kështu që, 
\[ \sum_{j=0}^n\ (-1)^j\s_jz^j\ \cdot \
\sum_{\nu=1}^\infty\ S_\nu\ z^\nu \ \ = \ \ \sum_{\mu=1}^n\ (-1)^{\mu+1}\mu\s_\mu z^\mu\]
 Duke krahasuar koefi\c cientët marrim atë,   që duam.

\qed

\begin{exa}
Çdo polinom kubik $f(x)$ mund të shkruhet në  formën 
\[ f(x) \ = \ x^3\ +\ a x\ + \ b\]
Vërtetoni se
\[ \D_f\ = \ -4a^3\ - \ 27 b^2.\]
 Kështu që,  formulat Kardano bëhen
\[ x_i\ = \ \left(-\frac{b}{2}\ +\ \sqrt{\frac{-\D_f}{108}} \right)^{1/3} \ +\ \left(-\frac{b}{2}\ -\ \sqrt{\frac{-\D_f}{108}} \right)^{1/3} 
\]
\end{exa}

\section{Format quadratike dhe një hyrje në algjebrën lineare} 

Teoria e formave është një nga më të vjetrat dhe më të bukurat e matematikës.  Ka vlera të pazëvendësueshme nga ana metodike sepse motivon përkufizimin e matricave (në fakt matricat lindën pikërisht nga format kuadratike), një pjesë të mirë të terminollogjise së algjebrës lineare (p.sh. definitisht pozitive, matricat simetrike, etj), dhe na jep ilustrime të shkëlqyera të aplikimit të algjebrës si per shembull klasifikimi i prerjeve konike, klasifikimi i sipërfaqeve algjebrike, etj.

Një  \textbf{formë binare kuadratike}  \es një polinom homogjen i gradës së dytë me dy ndryshore, pra një polinom i formës
\begin{equation}\label{form-2}
 f(x, y) = ax^2+ bxy + cx^2 
 \end{equation}
Pra \es thjesht polinomi kuadratik ku ne fusim një ndryshore të re $y$ dhe i bëjmë të gjitha termat me gradë totale dy.   Në fakt ky proçes \es i rëndësishëm në matematikën e lartë dhe quhet \textbf{homogjenizim} i polinomeve. Më poshtë ne do të spjegojmë këtë proçes për polinomet e gradës më të lartë, por për momentin le të përqëndrohemi tek polinomet kuadratike. 
Ne duam të studiojmë këto polinome kuadratike dhe rrënjet e tyre.

Le të vërejmë fillimisht disa veti të formave kuadratike. 

\begin{lem}
Për çdo dy forma kuadratike $f(x, y)$ dhe $g(x, y)$ shuma e tyre 
\[ (f+g ) (x, y) := f(x, y) + g(x, y) \]
\es përsëri një formë kuadratike. Për çdo konstante $\l \in \C$, polinomi $h(x, y):=\l f(x, y)$ \es një formë kuadratike. 
\end{lem}

Vërtetimi i kësaj leme \es elementar, por rëndësia e saj \es e madhe.  Ky rezultat elementar na jep përkufizimin e parë të një hapësire vektoriale.  Pra bashkësia e gjithë formave kuadratike me koefiçentë në $k$, ku $k$ \es secila prej $\Q, \R, \C$ \es një hapësirë vektoriale edhe pse një nxënës i klasës IX apo X nuk e ka dëgjuar më parë këtë koncept.

Një formë kuadratike $f(x, y)$ si në Eq.~\eqref{form-2} përcaktohet në mënyrë të vetme nga treshja e renditur e numrave $(a, b, c)$ dhe nga çifti i renditur i ndryshoreve $(x, y)$.  Në një farë mënyre ne duam që të organizojme këto treshe të renditura në se duam të studiojmë format kuadratike.  Gauss dhe me vonë Hermite filluan ti vinin këto treshe në tabela të tipit
\[ \begin{bmatrix} a & b \\ b & c \end{bmatrix}, \]
e cila quhet \textbf{matricë}.  
Në fakt Gauss filloi të përdorte    $M$ dhe  $\v$ si më poshtë
\[ M = \begin{bmatrix} a & \frac b 2 \\ \frac b 2 & c \end{bmatrix}, \quad \textit{ and } \quad \v=\cvec{ x \\ y } \]
Kjo terminollogji ishte mjaft efiçente me marrëveshjen që 
\[ \begin{bmatrix} a & b \\ b & c \end{bmatrix} \cvec{ x \\ y } = \cvec{ ax + by \\ cx + dy }\]
dhe 
\[ [x, y]   \begin{bmatrix} a & b \\ b & c \end{bmatrix}  =    [ax + by,   bx + cy ]\]
Në vend të $[x, y]$ ne shpesh përdorim simbolin $\v^t$ dhe e quajmë \textbf{transpose} të $\v$-së.  Atëherë forma kuadratike $f(x, y)$ jepet si më poshtë 
\[f (x, y) = \v^t M \v = (x, y) \,  \begin{bmatrix} a & \frac b 2 \\ \frac b 2 & c \end{bmatrix} \, \cvec{ x \\ y } \ = \ ax^2+ bxy + cx^2.   \]
Pra, kemi një korespondencë biunivoke midis formave kuadratike dhe matricave të formës 
\[ \begin{bmatrix} a & r \\ r & c \end{bmatrix} \]
Për një formë të dhënë $f(x, y)$ matrica koresponduese shënohet me $M_f$.  

Pozicionet në një matricë shënohen me $(i, j)$ ku $i$ tregon numrin e rradhës dhe $j$ numrin e kolonës.  Matricat e mësipërme $M$ quhen matrica $2 \times 2$ meqënëse kanë 2 rradhë dhe 2 kolona.  Një matricë $2\times 2$ quhet \textbf{matricë simetrike} kur termat në pozicionet $(1, 2)$ dhe $(2, 1)$ janë të barabarta.  Pra matricat tona që korespondojnë me format kuadratike janë matrica $2 \times 2$ dhe simetrike.

Perfundiminsht, kemi rezultatin si më poshtë:
\begin{lem}
Ekziston një korespondencë biunivoke midis formave kuadratike dhe matricave $2\times 2$, simetrike.  Më konkretisht, kjo korespondencë jepet si më poshtë
\[ f(x, y) = ax^2 + bxy + cy^2  \;  \mapsto  M_f = \begin{bmatrix} a & \frac b 2 \\ \frac b 2 & c \end{bmatrix}. \]
Matrica $M_f$ quhet \textit{matrica koresponduese} e formës $f(x, y)$. 
\end{lem}

Shihni pra se pa njohuri shtesë dhe pa shumë punë ne mund të prezantojmë nxënësin me konceptin e polinomeve homogjenë, matricave, matricave simetrike, shumëzimit të matricave, hapësirës vektoriale. Pra thjesht polinomi i gradës së dytë \es një minerë floriri.  Dhe ne vetëm sa kemi filluar.

Diskriminanti if një forme kuadratike përkufizohet njësoj si diskriminanti i polinomit $f(x, 1)$.  Pra, diskriminanti i $f(x, y)$ dhënë në Eq.~\eqref{form-2} \es $\D_f= b^2 - 4ac$.

Për një  matricë $A = \begin{bmatrix} a_{1, 1}  & a_{1, 1} \\ a_{2, 1}  & a_{2, 2} \end{bmatrix}$  ne përkufizojmë \textbf{determinantin} $\det A$  (ose \textbf{përcaktorin} siç përdoret në Shqip) si 
\[ \det A = a_{1, 1} a_{2, 2} - a_{2,1} a_{1, 2} \]
Atëherë, 
\[ \det M_f = ac - \frac {b^2} 4 \]
Vërtetimi i Lemës së mëposhtme tani \es një ushtrim elementar. 
\begin{lem}
Diskriminanti $\D_f$ i një forme kuadratike    $f(x, y) $ \es zero \att  $\det M_f =0$.  Për më tepër,
\[ \D_f = - 4 \, \det M_f\]
\end{lem}

Ka shumë autorë që i përkufizojnë format kuadratike si 
\[ f(x, y) = a x^2 + 2b \, xy + c y^2   \]
në mënyrë që matrica koresponduese \es  $M_f = \begin{bmatrix} a & b \\ b & c \end{bmatrix}$.  Atëherë,  diskriminanti  $\D_f = ac-b^2$ në vend të $b^2-4ac$.  
Kjo \es thjesht çështje preference, në se autori preferon ta nisë nga format kuadratike apo nga matricat.


\subsection{Ekuivalenca e formave, matricat e ngjashme}

I kthehemi edhe njëherë problemit të prerjeve konike, pra Prob.~(10).  Jepet një prerje konike
\[ ax^2 + b xy  + cy^2 + dx + ey = \l \]
Si duhet ndryshuar sistemi koordinativ , pra $\cvec{ x \\ y}$ që kjo prerje konike të jetë një nga format standart.  Së pari vemë re se grada e ekuacionit të ri nuk mund të ndryshojë pasi në këtë rast nuk do të kishim më një prerje konike. Pra trasformimet e mundshme janë vetëm transformimet 
\[ T:  (x, y)    \mapsto  \left (      \l_1 x + \l_2 y ,   \l_3 x + \l_4 y \right)  \]
Me fjalë të tjera, ne kemi një sistem të ri koordinatash që jepet nga 
\[
\cvec{ x^\prime \\ y^\prime } =  \begin{bmatrix} \l_1 & \l_2 \\ \l_3 & \l_4 \end{bmatrix} \cvec{ x \\ y }. 
\]
Funksioni $T(x, y)$ quhet \textbf{transformim linear}  dhe matrica $C= \begin{bmatrix} \l_1 & \l_2 \\ \l_3 & \l_4 \end{bmatrix}$ quhet \textbf{matrica e transformimit} $T$.  
$T(x, y)$ është funksion bijektiv dhe i anasjellti i tij $T^{-1}$ \es gjithashtu linear. Pra ekziston nje matricë për $T^{-1}$ qe ne e shënojmë me $C^{-1}$. Matrica $C^{-1}$ quhet matrica e anasjelltë e matricës $C$.  
Atëherë paraqisim edhe një herë problemin tonë:

\begin{prob}
Gjeni numrat $\l_1, \l_2, \l_3, \l_4$ të tillë që prerja konike transformohet në ekuacionin standart.
\end{prob}

Kjo motivon përkufizimin e mëposhtëm.  Dy forma kuadratike  $f(x, y)$ dhe  $g(x, y)$ do të quhen  \textbf{ekuivalente} \nqs ekzistojnë  $\l_1, \dots , \l_4 \in \C$ të tillë që 
\[ f  \left(\l_1 x + \l_2 y, \, \l_3 x + \l_4 y \right) = g(x, y).\]
Vini re se  
\[ 
\begin{split}
 f  \left(\l_1 x + \l_2 y, \, \l_3 x + \l_4 y \right) & = \x^t \, M_f \, \x  
= \left( \begin{bmatrix} \l_1 & \l_2 \\ \l_3 & \l_4 \end{bmatrix} \cdot \cvec{ x \\ y }\right)^t \, M_f \, \left( \begin{bmatrix} \l_1 & \l_2 \\ \l_3 & \l_4 \end{bmatrix} \cdot \cvec{ x \\ y } \right) \\
& = [x, y] \, \left( \begin{bmatrix} \l_1 & \l_2 \\ \l_3 & \l_4 \end{bmatrix}^t \, M_f \, \begin{bmatrix} \l_1 & \l_2 \\ \l_3 & \l_4 \end{bmatrix} \right) \,  \cvec{ x \\ y } \\
\end{split}
\]
Pra 
\[ f  \left(\l_1 x + \l_2 y, \, \l_3 x + \l_4 y \right) = \x^t \, M_f \, \x \]
ku 
\[ \x = \begin{bmatrix} \l_1 & \l_2 \\ \l_3 & \l_4 \end{bmatrix} \cdot \cvec{ x \\ y } \]
Matricat koresponduese të dy formave ekuivalente i quajmë \textbf{matrica të ngjashme}.  

\begin{exe}
Dy matrica $A$ dhe $B$ janë të ngjashme \att ekziston një matricë $C$ e tillë që 
\[ A = C^{-1} B C \]
\end{exe}

Ushtrimi i mëposhtëm \es standart në shkollat tona. 
\begin{exa}
Për një polinom kuadratik 
\[ f(x) = a x^2 + b x + c \]
me koefiçentë realë,  shenja e vlerës së $f(x)$ përcaktohet si vijon:
   $f(x) $ ka shenjën e kundert të $a$-së  në intervalin  $(-\alpha_1, \alpha_2)$ dhe ka shenjën e  $a$-së kudo tjetër.
\end{exa}
\begin{table}[htp]
\begin{center}
\begin{tabular}{c|c|ccc|c}
x    &     & $\alpha_1$   &           & $\alpha_2$ & \\
\hline
f(x)  & \qquad  a \qquad  &             & -a         &           &  \qquad     a \qquad  \\
\end{tabular}
\end{center}
\vspace{.5cm}
\caption{Studimi i shenjës së polinomeve kuadratike.}
\label{default}
\end{table}
Një formë kuadratike $f(x, y)$ quhet  \textbf{definitisht pozitive} \nqs  $f(x, 1) > 0$ për çdo $x\in \R$.  
\begin{exe}
$f(x, y)$ \es definitisht pozitive \att  $a> 0$ dhe $\D_f < 0$. 
\end{exe}

Për një matricë çfardo 
\[ A = \begin{bmatrix} a & b \\ c & d \end{bmatrix} \]
ne quajmë \textbf{eigenvlera} të matricës zgjidhjet e ekuacionit
\[ (\l - a)(\l-c) - bc =0 \]
%

\subsection{Klasifikimi i prerjeve konike}

Le të përpiqemi tani ti përgjigjemi Pyetjes~\ref{pyetje:conic}.  Pra jepet Equacioni~\eqref{conic}, ç'mund të themi për formën e grafikut? 

Në se ekuacioni në \eqref{conic} do të kishte $b=0$ \at ky do të ishte një ushtrim elementar.  Ne plotësonim katrorin për $ax^2+dx$ si edhe katrorin për $cy^2+ey$ dhe do të merrnim njue ekuacion të formës 
\[ A (x + \a)^2  + B (x+\b)^2 = C. \]
Me zëvendësimet $X=x+\a$ dhe $Y = y + \b$ ne kemi 
\[ A X^2 + B Y^2 = C,\]
i cili \es bashkësi boshe në $\R^2$ kur $C< 0$,  një elips kur $A, B$ janë me të njëjtën shenjë, dhe një hiperbolë kur $A, B$ janë me shenja të ndryshme.

Forma kuadratike
\[ G(x, y) = \l_1  x^2 + \l_2  y^2, \]
quhet \textbf{ formë diagonale} dhe matrica koresponduese 
\[ M_g = \begin{bmatrix} \l_1 & 0 \\ 0 & \l_2 \end{bmatrix} \]
\es njue matricë diagonale. 

Pra termi $b\, xy $ \es çfare ne duam të bejmë zero që të përcaktojmë formën e grafikut.  Ne mund të supozojmë që pas plotësimit të katrorëve dhe zëvendësimeve përkatëse, ekuacioni na  \es dhënë në formën  $ax^2+b xy + cy^2 = \l$.

Pra kemi problemin e mëposhtëm 

\begin{prob}
Jepet forma kuadratike 
\[ F(x, y) = A x^2+ B xy + C y^2. \]
Gjeni  zëvendësimet e nevojshme algjebrike (pra ndryshimin e sistemit koordinativ) 
\[ x = a x + b y, \qquad y=  c x + d y \]
që   forma $G(x, y) = F( a x + b y, c x + d y) $ \es diagonale. 
\end{prob}

Ne po i shmangim detajet e zgjidhjes së këtij problemi pasi ky problem \es thjesht një rast i vaçantë i problemit pasardhës. Po anallogjia \es e qartë, matrica $M_F$ diagonalizohet ne mënyrë ortogonale, pra 
\[ M_F = \begin{bmatrix} a & b \\ c & d \end{bmatrix}^{-1} \;  \begin{bmatrix} \l_1 & 0 \\ 0 & \l_2 \end{bmatrix} \; \begin{bmatrix} a & b \\ c & d \end{bmatrix} \]
e tillë qe $\l_1, \l_2$ janë eigenvlerat e $M_F$.   Lema në vijim përcakton formën e grafikut të  $F(x, y)=k$, për çdo konstant $k \in \R^\ast$. 

\begin{lem}
Grafiku 
\[ F(x, y) = ax^2+b xy + cy^2=k\]
\es elipse \nqs të dy egeinvalutat e $M_F$ janë pozitive dhe hiperbolë në se njëra \es pozitive dhe tjetra negative. 
\end{lem}
Ne i sygjerojmë lexuesit të shoh \cite{lin-alg} për detajet.


\subsection{Sipërfaqet kuadratike dhe klasifikimi i tyre}  
Ne mund të homogjenizojmë ekuacionin e dhënë në Ek.~\eqref{conic} si më poshtë 
\[ ax^2+ b y^2+ c z^2 + d xy + e xz + f yz \] 
duke futur një ndryshore të re $z$.  Pra ne kalojmë nga plani $\R^2$ në sistemin në hapësirë $\R^3$.  Prerjet tona konike tashmë janë thjesht projeksione në plan të grafikut  të sipërfaqes së mësipërme.   Pa  ndonjë kusht shtesë ne mund të supozojmë se  koefiçentët $d, e, f$ janë $2d, 2e, 2f$.

Një \textbf{formë ternare kuadratike} quhet polinomi homogjen kuadratik me tre ndryshore $x$, $y$, $z$. Pra një polinom i formës 
\[F(x,y,z)= ax^2+ b y^2+ c z^2 + 2 d xy + 2 e xz + 2 f yz \]
ku koefiçentët  $a, b, c, d, e, f$ janë numra realë.   Konsiderojmë ekuacionin 
 \[F(x, y, z)= h.\]
për ndonjë $h \in \R$. 
Në mënyrë plotësisht të ngjashme me format binare, ky ekuacion mund të shkruhet
\[  F(x, y, z) \; = \;    \x^t \, M_F \, \x  \]
ku
\[
\x =\cvec{ x \\ y \\ z  } \quad \textit{dhe} \quad
M_F = \left[
\begin{tabular}{rrr}
 a & d & e \\
 d & b& f\\
 e & f & c\\
\end{tabular}
\right]
\]
$M_F$ quhet matrica koresponduese e $F(x, y, z)$. Teoria e formave binare përgjithësohet në këtë rast fjalë për fjalë.  

\begin{pyetje} Ç'mund të themi për grafikun 
 \[F(x, y, z)= h,\]
në varësi të koefiçentëve të $F(x, y, z)$.
\end{pyetje}

Grafiku \es një sipërfaqe kuadratike në $\R^3$.  Një përshkrim i detajuar i këtyre sipërfaqeve jepet në \cite{kalk}.  Ne po i përkufizojmë më poshtë shkurtimisht. \\

%
\noindent \begin{minipage}{0.5\textwidth}
Një prej llojeve të sipërfaqeve kuadratike është    \textbf{elipsoidi}, i cili jepet me ekuacionin 
\[\frac{x^2}{a^2}+\frac{y^2}{b^2}+\frac{z^2}{c^2}=1\]
Në rastin kur $a^2=b^2=c^2$ ekuacioni i elipsoidit  paraqet një  sferë. Prerjet tërthore të tij me planet koordinative janë elipsa.
Elipsoidi është një nga sipërfaqet më të pergjithshme kuadratike, sferoidi dhe sferat janë raste të veçanta të ellipsoidit. 

\end{minipage}
\begin{minipage}{0.5\textwidth}
\flushright
\begin{tikzpicture}[scale=.9]
 \usetikzlibrary{arrows}
 \definecolor{ellipsecolor}{HTML}{AAAAFF}
 \shade [ball color=ellipsecolor] (0,0) ellipse (2.8 and  1.5);
 \draw [line width=0.2pt] (-2.8,0) arc (180:360:2.8 and  .35);
 \draw [dashed,line width=0.2pt] (2.8,0) arc (0:180:2.8 and  .35);
 \draw [line width=0.2pt] (0,1.5) arc (90:270:.4 and  1.5);
 \draw [dashed,line width=0.2pt] (0,-1.5) arc (-90:90:.4 and  1.5);
 \draw [black!60,line width=0.3pt,-latex] (0,0) -- (3.5,0,0);
 \draw [black!60,line width=0.3pt,-latex] (0,0) -- (0,2,0);
 \draw [black!60,line width=0.3pt,-latex] (0,0) -- (0,0,4.5);
 \pgfputat{\pgfpointxyz{3.1}{0.2}{0}}{\pgfbox[center,center]{\small y}};
 \pgfputat{\pgfpointxyz{0.2}{1.9}{0}}{\pgfbox[center,center]{\small z}};
 \pgfputat{\pgfpointxyz{0.2}{0}{4.3}}{\pgfbox[center,center]{\small x}};
 \pgfputat{\pgfpointxyz{0.05}{-0.2}{0}}{\pgfbox[center,center]{\small 0}};
 \node [below,right] at (-0.49,-0.49) {$a$};
 \node [below,right] at (2.7,-0.2) {$b$};
 \node [above,left] at (0,1.65) {$c$};
\end{tikzpicture}

\end{minipage}

\bigskip

Ekziston një tjetër elipsoid që quhet imagjinar dhe ka ekuacion
\[ \frac{x^2}{a^2}+\frac{y^2}{b^2}+\frac{z^2}{c^2}= - 1\]
Për detaje të mëtejshme shih \cite{kalk}.


\bigskip

\noindent \begin{minipage}{0.5\textwidth}

\centering
\begin{tikzpicture}
 \usetikzlibrary{arrows}
 \definecolor{insideo}{HTML}{798084}
 \definecolor{insidei}{HTML}{F0F0F0}
 \definecolor{outer}{HTML}{424296}
 \definecolor{inner}{HTML}{D8D8FF}
 \shadedraw [left color=insideo,right color=insideo,middle color=insidei] (0,3) ellipse (2 and  0.7);
 \shadedraw [left color=outer,right color=outer,middle color=inner]
 (2,3) arc (360:180:2 and  0.7) -- (-2,3) parabola bend (0,0) (2,3);
 \draw [dashed,line width=0.2pt] (-2,3) -- (2,3);
 \draw [dashed,line width=0.2pt] (-0.66,2.34) -- (0.66,3.66);
 \draw [line width=0.2pt] (-0.66,2.34) parabola [bend at end] (0,0);
 \draw [dashed,line width=0.2pt] (0.66,3.66) parabola [bend at end] (0,0);
 \draw [line width=0.2pt] (-1.25,1.2) arc (180:360:1.25 and  0.4);
 \draw [dashed,line width=0.2pt] (-1.25,1.2) arc (180:0:1.25 and  0.4);
 \draw [black!60,line width=0.3pt,-latex] (-2,0) -- (2,0,0);
 \draw [black!60,line width=0.3pt,-latex] (0,0) -- (0,4,0);
 \draw [black!60,line width=0.3pt,-latex] (0,0) -- (0,0,2);
 \pgfputat{\pgfpointxyz{1.9}{0.2}{0}}{\pgfbox[center,center]{\small y}};
 \pgfputat{\pgfpointxyz{0.2}{3.9}{0}}{\pgfbox[center,center]{\small z}};
 \pgfputat{\pgfpointxyz{0.2}{0}{1.8}}{\pgfbox[center,center]{\small x}};
 \pgfputat{\pgfpointxyz{0.05}{-0.2}{0}}{\pgfbox[center,center]{\small 0}};
\end{tikzpicture}

\end{minipage}
\begin{minipage}{0.5\textwidth}

\textbf{Paraboloidi eliptik} është    një  tjetër   sipërfaqeje kuadratike ekuacioni i së cilës është    i formës:
\[
\qquad \frac{x^2}{a^2}+\frac{y^2}{b^2}=\frac z c
\]

Prerjet me planet paralele me planin $xy$ janë elipse, ndërsa prerja me vetë planin $xy$ është një  pikë e vetme. 

Figura tregon rastin kur $c>0$, ndërsa  kur $c<0$, sipërfaqja është e kthyer me kokë poshtë. Në rastin kur $a=b$, sipërfaqja është një cilindër. \\

\end{minipage}
%

\bigskip

\noindent \begin{minipage}{0.5\textwidth}
Një tjetër sipërfaqe kuadratike   është    \textbf{paraboloidi hiperbolik} që   jepet me ekuacionin
\[
\frac{x^2}{a^2}-\frac{y^2}{b^2}=\frac z c
\]

Pra, një nga ndryshoret është e gradës së parë kurse pjesa tjetër është diferenca e dy katrorëve të ndryshoreve të tjera. Paraboloidi hiperbolik jep një shembull të atyre që ne i quajmë   \textbf{pika shalë}
 të cilat janë edhe maksimume lokale edhe minimume lokale; shih kapitujt mbi analizën me disa ndryshore në \cite{kalk}. 

\end{minipage}
\begin{minipage}{0.5\textwidth}
\centering
   \includegraphics[width=3in]{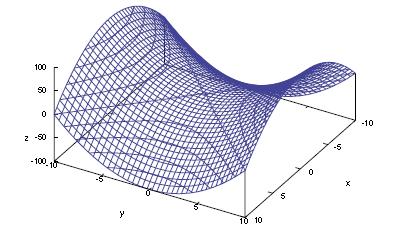} 
\end{minipage}

Ne sygjerojmë ushtrimin e mëposhtëm.

\begin{exe}
Diskutoni se çfarë vëndi gjeometrik  mund të jetë prerja e një sferoidi me një paraboloidi hiperbolic.  A është kjo gjithmonë një kurbë konike?
\end{exe}

%
\noindent \begin{minipage}{0.5\textwidth}

Një lloj tjetër  sipërfaqesh kuadratike \es \textbf{hiperboloidi me një  fletë} i cili jepet me ekuacionin:
\[
\frac{x^2}{a^2}+\frac{y^2}{b^2}-\frac{z^2}{c^2}=1
\]

Për hiperboloidin me një  fletë, prerja tërthore me \cdo plan paralel me planin $xy$ është    një  elips, ndërsa prerjet tërthore me plane paralele me planet $xz$ ose $yz$ janë hiperbola. Përjashtim bëjnë rastet e veçanta kur $x=\pm a$ dhe $y=\pm b$; në   këto plane prerjet janë çifte drejtëzash prerëse.

\end{minipage}
\begin{minipage}{0.5\textwidth}

\centering
  \begin{tikzpicture}
  \usetikzlibrary{arrows}
  \definecolor{insideo}{HTML}{798084}
  \definecolor{insidei}{HTML}{F0F0F0}
  \definecolor{outer}{HTML}{424296}
  \definecolor{inner}{HTML}{D8D8FF}
  \shade [left color=insideo,right color=insideo,middle color=insidei] (2,2) arc (0:180:2 and  .5) --
   (-2,2) arc (180:360:2 and  .5);
  \shadedraw [left color=outer,right color=outer,middle color=inner] (-2,-2) arc (180:360:2 and  .5).. controls
   (1,-1) and  (1,1).. (2,2)  -- (2,2) arc (360:180:2 and  .5).. controls (-1,1) and  (-1,-1).. (-2,-2);
  \draw (2,2) arc (0:180:2 and  .5);
  \draw [line width=0.2pt] (-1.25,0) arc (180:360:1.25 and  .4);
  \draw [dashed,line width=0.2pt] (1.25,0) arc (0:180:1.25 and  .4);
  \draw [dashed,line width=0.2pt] (2,-2) arc (0:180:2 and  .5);
  \draw [line width=0.2pt] (-0.8,-2.47).. controls (-0.21,-1.485) and  (-0.21,0.515).. (-0.8,1.53);
  \draw [dashed,line width=0.2pt] (0.8,-1.53).. controls (0.21,-0.515) and  (0.21,1.485).. (0.8,2.47);
  \draw [dashed,line width=0.2pt] (-0.8,-2.47) -- (0.8,-1.53);
  \draw [dashed,line width=0.2pt] (-0.8,1.53) -- (0.8,2.47);
  \draw [dashed,line width=0.2pt] (-2,2) -- (2,2);
  \draw [dashed,line width=0.2pt] (-2,-2) -- (2,-2);
  \draw [black!60,line width=0.3pt,-latex] (-2,0) -- (2,0,0);
  \draw [black!60,line width=0.3pt,-latex] (0,-3) -- (0,3,0);
  \draw [black!60,line width=0.3pt,-latex] (1,1) -- (0,0,7);
  \pgfputat{\pgfpointxyz{1.9}{0.2}{0}}{\pgfbox[center,center]{\small y}};
  \pgfputat{\pgfpointxyz{0.2}{2.9}{0}}{\pgfbox[center,center]{\small z}};
  \pgfputat{\pgfpointxyz{0.2}{0}{6.8}}{\pgfbox[center,center]{\small x}};
  \pgfputat{\pgfpointxyz{0.05}{-0.2}{0}}{\pgfbox[center,center]{\small 0}};
 \end{tikzpicture}
\end{minipage}

\textbf{Hiperboloidi me dy fletë}  është sipërfaqja kuadratike ekuacioni i të cilit \es 
\[ \frac{x^2}{a^2}-\frac{y^2}{b^2}-\frac{z^2}{c^2}=1\]
%


\noindent \begin{minipage}{0.5\textwidth}
\centering
  \begin{tikzpicture}[scale=.8]
  \usetikzlibrary{arrows}
  \definecolor{insideo}{HTML}{798084}
  \definecolor{insidei}{HTML}{F0F0F0}
  \definecolor{outer}{HTML}{424296}
  \definecolor{inner}{HTML}{D8D8FF}
  \shadedraw [rotate=-45,left color=insideo,right color=insideo,middle color=insidei,shading angle=45]  (0,-2.5) ellipse (1 and  0.5);
  \shadedraw [rotate=-45,left color=outer,right color=outer,middle color=inner,shading angle=45]
  (-1,-2.5) arc (180:0:1 and  0.5) -- (1,-2.5) parabola bend (0,-0.5) (-1,-2.5);
  \shadedraw [rotate=-45,left color=outer,right color=outer,middle color=inner,shading angle=45]
  (-1,2.5) arc (180:0:1 and  0.5) -- (1,2.5) parabola bend (0,0.5) (-1,2.5);
  \draw [rotate=-45,dashed,line width=0.2pt] (-1,2.5) arc (180:360:1 and  0.5);
  \draw [rotate=-45,dashed,line width=0.2pt] (-0.4472,-2.9472) -- (0.4472,-2.0528);
  \draw [rotate=-45,dashed,line width=0.2pt] (0.4472,-2.9472) -- (-0.4472,-2.0528);
  \draw [rotate=-45,line width=0.2pt] (0.4472,-2.0528) parabola [bend at end] (0,-0.5);
  \draw [rotate=-45,dashed,line width=0.2pt] (-0.4472,-2.9472) parabola [bend at end] (0,-0.5);
  \draw [rotate=-45,dashed,line width=0.2pt] (-0.4472,2.9472) -- (0.4472,2.0528);
  \draw [rotate=-45,dashed,line width=0.2pt] (0.4472,2.9472) -- (-0.4472,2.0528);
  \draw [rotate=-45,dashed,line width=0.2pt] (-0.4472,2.0528) parabola [bend at end] (0,0.5);
  \draw [rotate=-45,line width=0.2pt] (0.4472,2.9472) parabola [bend at end] (0,0.5);
  \draw [black!60,line width=0.3pt,-latex] (-2.5,0) -- (2.5,0,0);
  \draw [black!60,line width=0.3pt,-latex] (0,-3) -- (0,3,0);
  \draw [black!60,line width=0.3pt,-latex] (2.5,2.5) -- (0,0,7);
  \pgfputat{\pgfpointxyz{2.4}{0.2}{0}}{\pgfbox[center,center]{\small y}};
  \pgfputat{\pgfpointxyz{0.2}{2.6}{0}}{\pgfbox[center,center]{\small z}};
  \pgfputat{\pgfpointxyz{0.2}{0}{6.8}}{\pgfbox[center,center]{\small x}};
  \pgfputat{\pgfpointxyz{0.05}{-0.2}{0}}{\pgfbox[center,center]{\small 0}};
  \end{tikzpicture}
\end{minipage}
\begin{minipage}{0.5\textwidth}
Për hiperboloidin me dy fletë, prerja tërthore me \cdo plan paralel me planet $xy$ ose $yz$ është    një  hiperbolë. Me planin $yz$ nuk ka prerje tërthore,  sepse për $x=0$ ekuacioni
\[  -\frac{y^2}{b^2}-\frac{z^2}{c^2}=1 \]
nuk ka zgjidhje.  Me \cdo plan paralel me planin $yz$ për  të cilin $|x|>a$, prerja është    elips.

Në vijim ne do të mësojmë se si të klasifikojmë gjithë sipërfaqet kuadratike sipas këtyre llojeve.  \\

\end{minipage}
 

\noindent \begin{minipage}{0.5\textwidth}

\textbf{Cilindri parabolik} \es sipërfaqja kuadratike ku një nga ndryshoret \es në fuqi të dytë dhe tjetra në fuqi të parë, ndryshorja e tretë nuk ekziston. Pra një sipërfaqe e tillë
 \[ x^2 = \lambda z\]

\end{minipage}
\begin{minipage}{0.5\textwidth}
   \centering
   \includegraphics[width=2in]{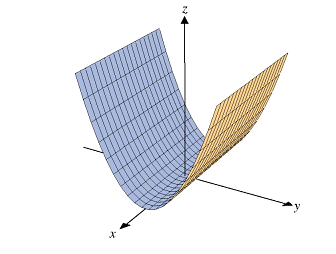} 
\end{minipage}

Megjithëse ka forma të tjera të sipërfaqeve kuadratike, secila prej tyre mund të përftohet si një degjenerim i formave të mësipërme. 

I kthehemi tani problemit tonë ku na jepet një sipërfaqe kuadratike e përgjithshme dhe duam të gjejmë llojin a saj në varësi të koefiçentëve. 
Si në rastin e prerjeve konike, \es e përshtatshme që të rrotullojmë sistemin koordinativ që termat 
 $xy, yz, xz$ të zhduken. Forma të tilla, pa termat $xy, yz, xz$,  quhen \textbf{forma kuadratike diagonale}. 



Përkufizojmë  \textbf{inercinë}    $\mbox{in } M_F$ e një forme quadratike $F(x, y, z)$  si treshen e renditur  
\[\mbox{in } M_F:=( n_1, n_2, n_3 ),\] 

ku  $n_i$,  $i=1,2,3 $ \es përkatësisht numri pozitiv, negativ, dhe zero i eigenvlerave të $M_F$-së.  Atëherë kemi lemën e mëposhtme.

\begin{lem}
Jepet forma kuadratike ternare  $F(x, y, z)$ dhe matrica koresponduese  $A$. Pohimet e mëposhtme janë të vërteta:
\begin{itemize}
\item[i)] \Nqs $\mbox{in } M_F= (3, 0, 0)$, atëherë sipërfaqja kuadratike është një  ellipsoid.

\item[ii)] \Nqs $\mbox{in } M_F= (2, 0, 1)$, atëherë sipërfaqja kuadratike është një   paraboloid eliptik.

\item[iii)] \Nqs $\mbox{in } M_F= (2, 1, 0)$, atëherë sipërfaqja kuadratike është një hyperboloid me një fletë.

\item[iv)] \Nqs $\mbox{in } M_F= (1, 2, 0)$, atëherë sipërfaqja kuadratike është një hyperboloid me dy fletë.

\item[v)] \Nqs $\mbox{in } M_F= (1, 1, 1)$, atëherë sipërfaqja kuadratike është një  paraboloid hiperbolik.

\item[vi)] \Nqs $\mbox{in } M_F= (1, 0, 2)$, atëherë sipërfaqja kuadratike është një cylinder parabolic.\\
\end{itemize}
\end{lem}    

Vërtetimi i Lemës mund të gjendet në \cite{lin-alg} dhe bazohet në faktin se numri i engenvlerave pozitive, negative, ose zero nuk ndryshon pavarësisht sistemit koordinativ të zgjedhur.  Në fund të fundit një konfirmim që forma e sipërfaqes nuk ndryshon, pra një elipsoid nuk mund të bëhet parabolois apo anasjelltas.  Një koncept që nxënësi e beson shumë më kollaj sepse bazohet në intuitën gjeometrike.

Më poshtë po japim shkurtimisht metodën se si përkaktohet forma diagonale dhe zevendësimet algjebrike që e bëjnë këtë të mundur.

\begin{prob}
Gjeni diagonalizimin ortogonal të një matrice simetrike $A$. 
\end{prob}

Pra, ne duam të gjejmë një matricë ortogonale    $S$ dhe një matricë diagonale $D$ të tillë që  
\[ A= S^T D S.\] 
Kujtoni që për matricat ortogonale  $S^T= S^{-1}$; shihni \cite{lin-alg} për detajet. 

Së pari gjejmë gjithë eigenvlerat e $A$-së 
\[ \lambda_1, \dots , \lambda_r, \]
dhe shumfishmëritë e tyre.

Së dyti, për çdo eigenvlerë $\lambda_i$ gjejmë një bazë ortonormale 
\[ \mathcal B = \{ v_{i, 1}, \dots , v_{i, s_i} \}  \]
Së treti, krijojmë matricën 
\[ S = \left[ v_{1,1} \, | \, \dots \, | v_{1, s_1} | v_{2, 1} | \dots | v_{2, s_2} | \dots | v_{r, s_r}   \right] \]
dhe matricën diagonale 
\[ D= \mbox{diag } ( \l_1, \dots , \l_r)   \]
që janë matricat e kërkuara.

Problemi i mësipërm ka vlerë për çdo matricë simetrike  (pra në mënyrë ekuivalente për çdo formë kuadratike, jo vetëm format binare dhe ternare).  Një nxënës që kupton gjithë  hapat e zgjidhjes së këtij problemi mund të thuhet se ka një bazë të shendoshë të algjebrës lineare. 

I përshtatur për sipërfaqet kuadratike ky problem bëhet: 

\begin{prob}Konsideroni formën kuadratike 
\[ q (x, y, z) = x^2+y^2+z^2+2xy+2xz+2yz\]

i) Gjeni matriocën koresponduese  $A$ të $q(x, y, z)$.  

ii) Gjeni matricat  $C$ dhe $D$ të tilla që  $A=C^{-1} D C$, ku $C$ \es ortogonale dhe  $D$ matricë  diagonale.

iii) Përcaktoni zëvendësimet lineare  
\[ 
\begin{split}
x^\prime & = a_1 x + b_1 y  + c_1 z , \\
  y^\prime & =a_2 x + b_2 y  + c_2 z , \\
   z^\prime & =a_3 x + b_3 y  + c_3 z , \\
\end{split}
\]
të tilla që  $q (x^\prime, y^\prime, z^\prime) $  transformohet në një formë diagonale që i korespondon  $D$-së. 

iv) Gjeni  $q (x^\prime, y^\prime, z^\prime) $ algjebrikisht për të kontrolluar që vërtet i korespondon matricës diagonale $D$.

v) Çfarë \es forma e sipërfaqes   $q (x^\prime, y^\prime, z^\prime)  =4$?  Ç'mund të thoni për formën e sipërfaqes   $q (x, y, z) =4$? 
\end{prob}

Kur ushtrimi i mësipërm \es metodollogjik dhe përmban brenda të paktën një simestër të algjebrës lineare dhe një simestër të gjeometrisë analitike, të gjitha keto njohuri mund të bëhen në shkollën e mesme dhe nxënësi arrin nivelin e një studenti të vitit të dytë të universitetit. 

Problemi i mëposhtëm, në dukje më inoçent, \es një shembull konkret i metodës së mësipërme.

\begin{prob}
Klasifikoni sipërfaqen kuadratike
\[x^2+ y^2 - z^2 + 3xy- 5 xz + 4yz=1.\]

\end{prob}

\section{Ndërtimet gjeometrike}

Këtu po japim një permbledhje elementare të materialit të trajtuar ne Kapitullin~14 \te \cite{alg}. Ne supozojmë njohuri elementare të algjebrës dhe konceptin e fushës. Në rast se ky koncept nuk është zhvilluar, mësuesi mund të zëvendësojë fjalën \text{fushë} me një nga bashkësitë $\Q, \R$, ose $\C$. 

Gjithashtu ne supozojmë se nxënësi kupton konceptin e zgjerimit të fushave ose nënfushe e një fushe të dhënë. Pra ne do të përdorim shtrirjet algjebrike të fushave.  Një material pregatitor për të lexuar këtë material \es materiali elementar në Kapitullin~13 \ne \cite{alg}. 
 
\subsection{Numrat e ndërtueshëm}

Një numër real $\a$  është  \textbf{i ndërtueshëm\/}\index{Constructible
number},   \nqs mund të ndërtojmë një segment me gjatësi $| \a |$ në një numër të fundëm
hapash nga një segment njësi,  duke përdorur një vizore dhe një kompas.

\begin{figure}[hb]
\begin{center}
\centerline {
\includegraphics[width=2in]{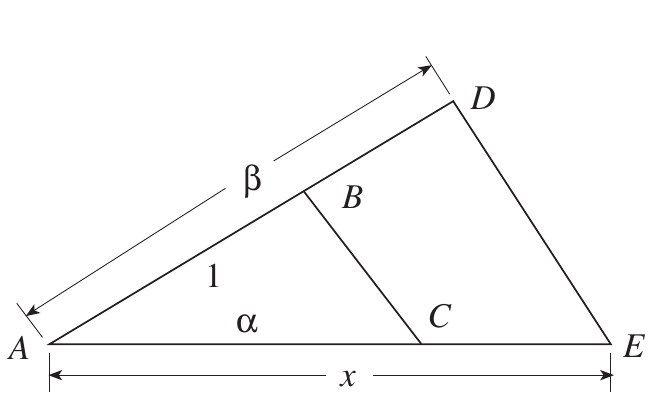}
}
\end{center}
\caption{Ndërtimi i  produkteve të numrave}
\label{Multiply}
\end{figure}

\begin{thm}
Bashkësia të gjithë numrave realë të ndërtueshëm formon një nënfushë $F$,   të fushës së numrave realë $\R$.
\end{thm}

\proof
Le të jenë $\a$ dhe $\beta$,  numra të ndërtueshëm. ne duhet të vërtetojmë se, 
$\a + \beta$,  $\a - \beta$,  $\a \beta$,  dhe $\a /
\beta$ ($\beta \neq 0$),   janë gjithashtu numra të ndërtueshëm. Ne mund të supozojmë,   që të dy
 $\a$ dhe $\beta$,  janë pozitivë,  ku $\a > \beta$. Është pothuajse e qartë
,  sesi të ndërtojmë $\a + \beta$ dhe $\a -
\beta$. Për të gjetur një segment me gjatësi $\a \beta$,  supozojmë,  që
 $\beta > 1$ dhe ndërtojmë trekëndëshin në Figurën~\ref{Multiply}, 
 i tillë që,  trekëndëshat $\triangle ABC$ dhe $\triangle ADE$ janë të ngjashëm.
Meqënëse $\a / 1 = x / \beta$,  segmenti $x$ ka gjatësi
$\a \beta$. Një  ndërtim i ngjashëm mund të bëhet,   \nqs $\beta <1$. Po  e lëmë
si ushtrim të vërtetoni se,  i njëjti trekëndësh mund të përdoret për të ndërtuar
 $\a / \beta$,  për $\beta \neq 0$.
\qed

\begin{lem}
\Nqs $\a$  është një numër i ndërtueshëm,  atëherë  edhe $\sqrt{\a}$  është
numër i ndërtueshëm.
\end{lem}

\begin{figure}[htb]
\begin{center}
\centerline {
\includegraphics[width=2in]{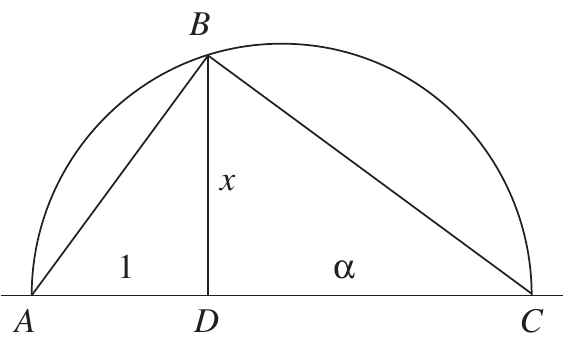}
}
\end{center}
\caption{Ndërtimi i rrënjëve \te numrave}
\label{Root}
\end{figure}

\proof
Në Figurën~\ref{Root} trekëndëshat $\triangle ABD$,  $\triangle BCD$, 
dhe $\triangle ABC$ janë të ngjashëm; pra,  $1 /x = x / \a$,  ose $x^2 =
\a$.
\qed

Lema e mëposhtme mund të jetë një detyrë e përshtatshme.  

\begin{lem}
Është i mundur ndërtimi i trekëndëshave të ngjashëm.
\end{lem}

Ne mund të gjejmë  në  plan \cdo pikë $P =( p,  q)$, 
 e cila ka koordinata racionale $p$ dhe $q$.  Ne duam të dimë se,   cilat pika të tjera
mund të ndërtohen me kompas dhe  vizore,   nga pikat me koordinata racionale.

\begin{lem}
Le të jetë $F$ një nënfushë e $\R$.
 
i)  \Nqs një drejtëz përmban dy pika  në $F$,  atëherë ajo ka ekuacionin $a x + by + c = 0$,  ku $a$,  $b$,  dhe $c$ janë në $F$.

ii)  \Nqs një rreth e ka qendrën në një pikë me koordinata në $F$ dhe me rreze,  cila  është gjithashtu në $F$,  atëherë ai ka ekuacionin $x^2 + y^2 + d x + e y + f = 0$,  ku $d$,  $e$,  dhe $f$ janë në $F$.
\end{lem}

\proof  Le të jenë $(x_1,  y_1)$ dhe $(x_2,  y_2)$,   pika të një drejtëze me koordinata
në $F$.  \Nqs $x_1 = x_2$,  atëherë  ekuacioni i drejtëzës,   që kalon nga dy pikat
  është  $x - x_1 = 0$,  i cili ka formën $a x + by +
c = 0$.  \Nqs $x_1 \neq x_2$,  atëherë  ekuacioni i drejtëzës,  që kalon nga dy pikat jepet nga
\[
y - y_1 = \left( \frac{y_2 - y_1}{x_2 - x_1} \right) (x - x_1), 
\]
i cili mund të kthehet në formën e duhur.

Për të vërtetuar pjesën e dytë të lemës,  supozojmë se,   $(x_1,  y_1)$  është  qendra
e një rrethi me rreze $r$.  Atëherë rrethi ka ekuacionin
\[
(x - x_1)^2 + (y - y_1)^2 - r^2 = 0.
\]
Ky ekuacion mund të shkruhet lehtë në formën e duhur.

\qed

Duke filluar me fushën e numrave të ndërtueshëm $F$,  kemi tre mënyra të ndryshme për të ndërtuar pikat shtesë  në $\R$ me  kompas dhe  vizore .
 
\begin{enumerate}

\item Për të gjetur pika e reja të mundshme në $\R$,  ne mund të marrim
prerjen e dy drejtëzave,  ku seicila prej tyre kalon nga
dy pika të dhëna me koordinata në $F$.

\item Prerja e një drejtëze,   e cila kalon nga dy pika me koordinata   në $F$ me një rreth, 
qendra e të cilit i ka koordinatat në $F$
me gjatësi rrezje në $F$,  do të na japë pika të reja  në $\R$.

\item Ne mund të përftojmë pika të reja në
 $\R$,  duke prerë dy rrathë,  qendrat e të cilëve
 i kanë koordinatat
 në $F$ dhe gjatësitë e rrezeve të tyre janë
 në $F$.

\end{enumerate}

Rasti i parë,  nuk na jep pika të reja në $\R$,  meqë zgjidhja e dy ekuacioneve
të formës $a x + by + c = 0$,   me koefi\c cientë në $F$,   gjithmonë do të jetë në $F$. Rasti i tretë mund të sillet
te rasti i dytë.
Le të jenë

\begin{eqnarray*}
x^2 + y^2 + d_1 x +e_1 x + f_1 = 0 \\
x^2 + y^2 + d_2 x +e_2 x + f_2 = 0
\end{eqnarray*}
ekuacionet e dy rrathëve,  ku $d_i$,  $e_i$,  dhe $f_i$ janë në $F$ për $i = 1,  2$. Këta rrathë kanë të njëjtën  prerje si, 
rrethi
$$
x^2 + y^2 + d_1 x +e_1 x + f_1 = 0
$$ dhe drejtëza
$$
(d_1 - d_2) x + b(e_2 - e_1)y + (f_2 - f_1) = 0.
$$
Ekuacioni i fundit,   është   ai i kordës,   që kalon në pikat e prerjes së dy rrathëve. Pra,   prerja e dy
rrathëve mund të sillet në rastin e prerjes së një drejtëze me një rreth.

Konsiderojmë rastin e  prerjes së një drejtëze me një rreth, 
duhet të përcaktojmë natyrën e zgjidhjeve të ekuacioneve.
\begin{eqnarray*}
a x + by + c & = & 0 \\
x^2 + y^2 + d x + e y + f & = & 0.
\end{eqnarray*}
\Nqs eleminojmë $y$ nga këto ekuacione,  përftojmë një ekuacion të formës
 $Ax^2 + B x + C = 0$,  ku $A$,  $B$,  dhe $C$ janë në $F$. Koordinata
$x$ e pikave të prerjes  jepet nga
$$
x = \frac{- B \pm \sqrt{B^2 - 4 A C} }{2 A}
$$
dhe  është  në $F( \sqrt{\a}\,  )$,  ku $\a = B^2 - 4 A C > 0$.
Kemi vërtetuar lemën në vazhdim.

\begin{lem}
Le të jetë $F$  fusha e numrave të ndërtueshëm. Atëherë pikat
e prerjes së drejtëzave  dhe rrathëve në $F$,   ndodhen në
fushën $F( \sqrt{\a}\,  )$ për ndonjë $\a$ në $F$.
\end{lem}

\begin{thm}
Një numër real $\a$  është  numër i ndërtueshëm,   atëherë dhe vetëm atëherë,   kur  ekziston një  varg fushash
$$
\Q = F_0 \subset F_1 \subset \cdots \subset F_k
$$
 të  tilla që,  $F_i = F_{i-1}( \sqrt{ \a_i}\,  )$ ku $\a \in F_k$.
Në ve\c canti,   ekziston një numër i plotë $k > 0$,   i tillë që,  $[\Q(\a) : \Q ] = 2^k$.
\end{thm}

\proof
Ekzistenca e $F_i$-ve dhe e $\a_i$-ve  është pasojë direkte e
 Lemës \se mësipërme dhe e faktit,  që
$$
[F_k: \Q] = [F_k : F_{k-1}][F_{k-1} : F_{k-2}] \cdots [F_1:
\Q ] = 2^k.
$$

\qed

\begin{cor}
Fusha e gjithë numrave të ndërtueshëm  është  një shtrirje algjebrike e
$\Q$.
\end{cor}

Si\c c mund ta shihni nga fusha e numrave të ndërtueshëm, 
jo   \cdo shtrirje algjebrike e një fushe  është shtrirje e fundme.

\subsection{Dyfishimi i një kubi dhe katrori i rethit}

Tani jemi gati të shqyrtojmë problemet klasike të dyfishimit
të kubit dhe  të kthimit të rrethit në katror.
Ne mund të përdorim fushën e numrave të ndërtueshëm për të
treguar me saktësi se,  kur një ndërtim i ve\c cantë algjebrik
mund të realizohet.

{\em Dyfishimi i kubit  është  i pamundur\index{Doubling the cube}}. Kur jepet
brinja e kubit,   është  e pamundur të ndërtosh me vizore dhe kompas
brinjën e kubit,  i cili ka dyfishin e vëllimit të kubit origjinal.
 Le të jetë kubi origjinal me brinjë me gjatësi 1 dhe me vëllim po 1. \Nqs mund të ndërtojmë një kub me vëllim
 2,  atëherë  ky kub i ri do ta ketë brinjën me gjatësi
$\sqrt[3]{2}$. Megjithatë,  $\sqrt[3]{2}$  është një rrënjë e polinomit të pathjeshtuar
 $x^3 -2$,  mbi $\Q$; pra, 
$$
[\Q(\sqrt[3]{2}\,  ) : \Q] =3
$$
Kjo është  e pamundur,  sepse 3 nuk  është  fuqi e 2.

\begin{thm}  
Është e pamundur,   që të dyfishojmë kubin.
\end{thm}

\proof Marrim një kub me vëllim 1. Për të dyfishuar kubin  duhet të ndërtojmë një  $x$,   të tillë që,   $$x^3=2.$$ Polinomi $$f(x)=x^3-2$$ është i pathjeshtueshëm mbi $\Q$ dhe prandaj,  për \cdo rrënjë $\a$
të $f(x)$,   kemi $[\Q(\a) : \Q]=3$,  e cila nuk është fuqi e 2.

\qed

\subsection{Kthimi i rrethit në katror} \index{Squaunazë  the circle}  
Supozojmë,   që kemi një rreth me rreze 1.  Sipërfaqja e rrethit është
 $\pi$; prandaj,  ne duhet  të ndërtojmë një katror me brinjë$\sqrt{\pi}$. Kjo  është  e pamundur,  sepse si $\pi$ dhe
$\sqrt{\pi}$ janë të dy transhendentë. Kështu që,  me përdorimin e  vizores dhe kompasit,  është e pamundur të ndërtosh një katror me të njëjtën sipërfaqe sa të rrethit.

\begin{thm} Është  e pa mundur ta kthesh rrethin në katror.
\end{thm}

\proof Jepen $r$ i ndërtueshëm dhe një rreth me rreze $r$. Duam të ndërtojmë një katror me brinjë
 $x$,  i tillë që,  $$x^2=\pi r^2$$
Meqënëse  $\pi$  nuk është as numër algjebrik,  atëherë rrënjët e ekuacionit të mësipërm nuk janë as
algjebrikë dhe prandaj nuk mund të jenë të ndërtueshme.

\qed

\subsection{Ndarja në tresh e një këndi}\index{Trisection of an angle}
Ndarja në tresh e një këndi të \c cfarëdoshëm  është  e pamundur.  Do të vërtetojmë se, 
  është  e pamundur të ndërtosh një kënd  $20^\circ$.  Për pasojë,  një kënd
$60^{\circ}$ nuk mund të ndahet në tresh. Si fillim,  duhet të llogaritim formulën e
 këndit trefish për kosinusin:
\begin{equation}\label{cos3a}
\cos 3 \theta  = 4 \cos^3 \theta - 3 \cos \theta.
\end{equation}
Këndi $\theta$ mund të ndërtohet,   atëherë dhe vetëm atëherë,  kur  $\a = \cos
\theta$  është   i ndërtueshëm. Le të jetë $\theta = 20^{\circ}$. Atëherë $\cos 3
\theta =  \cos 60^\circ = 1/2$. Nga formula e këndit të trefishtë për
kosinusin, 
$$
4 \a^3 - 3 \a = \frac{1}{2}.
$$
Kështu që,   $\a$  është një rrënjë e $8 x^3 - 6 x -1$. Ky polinom
nuk ka faktorë në $\Z[x]$,  pra   është i pathjeshtueshëm mbi $\Q[x]$. Kështu që,   $[\Q( \a ) : {\mathbb Q }] = 3$. Për pasojë, 
$\a$ nuk mund të jetë një numër i ndërtueshëm.

\begin{thm} Është e pamundur të ndahet një kënd në tresh.
\end{thm}

\proof Të ndash në tresh një kënd $3 \a$,   është njësoj si të ndërtosh $\cos \a$,  kur jepet $\cos 3\a$.
Ekuacioni \eqref{cos3a}
 na jep polinomin $$f(x)= 4 x^3 - 3 x - \cos 3 \a, $$ për  të cilin $\cos \a$ është rrënjë. Për disa vlera të $\a$ polinomi $f(x)$ është i pathjeshtueshëm dhe $[\Q( \cos \a) : \Q]=3$,  e cila nuk është
fuqi e  2.
Marrim $\a= 20^\circ$. Atëherë, 
$$f(x)= 8 x^3 - 6x -1 $$ është i pathjeshtueshëm mbi $\Q$. Pra,   në qoftë se,  $\a$ është rrënjë e $f(x)$,  atëherë $[\Q(\a):\Q]=3$,   e cila nuk është fuqi e 2.

\qed

\begin{exa}Përcakto në se këta kënde mund të ndahet në tresh. 

  i) Këndi $\b$,  i tillë që,  $\cos \b = \frac 1 3$.

  ii) $\b = 120^\circ$
\end{exa}

\sol 
Të ndash në tresh një kënd $\b=3 \a$,  është njësoj si të ndërtosh $\cos \a$,  kur jepet $\cos 3\a$. Ekuacioni
$$\cos 3\a = 4 \cos^3 \a - 3 \cos \a$$
na jep polinomin 
\[ f(x)= 4 x^3 - 3 x - \cos 3 \a,  \] për të cilin $\cos \a$ është rrënjë. \\

 i) Në qoftë se,  $\cos 3\a = \frac 1 3$,  atëherë $\cos \a$ është një rrënjë e polinomit
$$f(x)=4 x^3-3 x- \frac 1 3 $$
i cili është i pathjeshtueshëm. Atëherë,  $[\Q (\cos \a) : \Q]=3$,  i cili nuk  është fuqi e 2. Pra,   ky kënd
nuk mund të ndahet në tresh.

 ii) Në qoftë se,  $\b=120^\circ$,  atëherë $\cos \b= - \frac 1 2$. Atëherë $$f(x) = 4 x^3-3 x + \frac 1 2 $$
i cili është i pathjeshtueshëm mbi $\Q$ dhe,  si më lart,  këndi nuk mund të ndahet në tresh. 

\qed

\subsection{Ndërtimi i një shumëkëndëshi të rregullt}

\begin{thm}
 $n$-këndëshi është i ndërtueshëm,  atëherë dhe vetëm atëherë,  kur  $$n= 2^k \cdot p_1 \cdots p_s, $$ ku $p_i$ janë numrat e thjeshtë Fermat,  të dalluar, pra të formës $p= 2^{2^r} +1$.
\end{thm}

\proof Ndërtimi i një $n$-këndëshi është ekuivalent me ndërtimin e $\cos \frac {2\pi } n$. 
Shënojmë me 
$e_n = \cos \frac {2\pi } n + i \sin \frac {2\pi } n$ rrënjën primitive të njësisë. 
Atëherë, 
$\cos \frac {2\pi } n =  (\e_n + \e_n^{-1})/2$. 
Pra,   $\Q(\e_n)$ është një shtrirje e  $\Q( \frac {2 \pi} n)$.
 Ne e mbarojmë vërtetimin vetëm  për $n$ numër të thjeshtë $p$,  pjesa tjetër do të vërtetohet në kapitullin e fushës ciklotomike \ne \cite{alg}. 
 Pra,  $\cos \frac {2\pi} n$ është i ndërtueshëm,   në qoftë se,  $\frac {p-1} 2 = 2^r$ për ndonjë $r
\geq 0$. Kështu që,    kjo është e mundur vetëm  për numra të thjeshtë $p$,   të formës $p= 2^k +1$. Këta janë saktësisht numrat e thjeshtë Fermat dhe ata janë të formës $p= 2^{2^r} +1$. 
\qed

Disa problema të përshtatshëm për gjeometrinë e shkollës së mesme janë:

\begin{prob}
Vërtetoni se,   9-këndëshi i rregullt nuk  është  i ndërtueshëm me vizore dhe kompas  kurse 20-këndëshi i rregullt  është i ndërtueshëm.
\end{prob}

\begin{prob}
Vërtetoni se,   kosinusi i një grade ($\cos 1^\circ$)  është algjebrik mbi
$\Q$,   por nuk është i ndërtueshëm.
\end{prob}

\begin{prob}
Vërtetoni se   \nqs $\a$ dhe $\beta$ janë numra të ndërtueshëm, \te  tillë që   $\beta \neq 0$,  atëherë   i tillë është edhe $\a / \beta$.
\end{prob}

%

\begin{prob}  
Jepni një mënyrë gjeometrike për të ndërtuar një $n$-gon të rregullt,  për 
\[n=3,  4,  5,  6,  8,  10,  12,  15,  16,  17,  20,  24\]
\end{prob}

Ndërtimi i $n$-këndëshit të rregullt për $n=3, \dots , 6$ zakonisht \es trajtuar bë programet tona \te gjeometrisë.  Natyrisht, nga kjo rastet $n= 8, 10, 12, 16, 20, 24$ janë rrjedhime elementare.  Rasti $n=15$ \es një problem interesant për nxënësin e vitit të dytë të gjimnazit. Rasti $n=17$ \es më i vështiri. 


\nocite{*}
\bibliographystyle{plain}

\bibliography{2016-5}{}

\end{document}